\documentclass[12pt]{amsart}
\usepackage{amssymb}
\usepackage{mathrsfs}
\usepackage{graphicx}
\renewcommand\a{\alpha}
\renewcommand\b{\beta}

\renewcommand\d{\delta}
\newcommand\la{\lambda}

\newcommand\e{\eta}

\newcommand\s{\sigma}

\newcommand\vf{\varphi}

\newcommand\D{\Delta}

\newcommand\vL{\varLambda}

\newcommand{\ZZ}{\mathbb Z}
\newcommand{\NN}{\mathbb N}

\newcommand\BQ{\mathbf Q}

\newcommand\Bp{\mathbf p}
\newcommand\Bm{\mathbf m}

\newcommand\CA{\mathcal{A}}
\newcommand\CB{\mathcal{B}}
\newcommand\ZC{\mathcal{C}}

\newcommand\CI{\mathcal{I}}

\newcommand\CT{ \mathcal{T}}

\newcommand\FS{\mathfrak S}

\newcommand\Fs{\mathfrak s}

\newcommand\Ft{\mathfrak t}

\newcommand\wh{\widehat}
\newcommand\wt{\widetilde}

\newcommand\ol{\overline}

\newcommand\trreq{\trianglerighteq}

\newcommand{\lan}{\langle}
\newcommand{\ran}{\rangle}

\newcommand{\trr}{\triangleright }
\newcommand{\trl}{\triangleleft }
\newcommand{\ra}{\rightarrow }
\newcommand{\ot}{\otimes}

\newcommand\Hom{\operatorname{Hom}}
\newcommand\End{\operatorname{End}}

\newcommand\res{\operatorname{res}}

\renewcommand\Im{\operatorname{Im}}

\newcommand\Std{\operatorname{Std}}
\newcommand\Ext{\operatorname{Ext}}

\newcommand{\rad}{\operatorname{rad}}

\newcommand{\isom}{\,\raise2pt\hbox{$\underrightarrow{\sim}$}\,}
\newcounter{ichi}
\newcommand{\roi}{\roman{ichi}}
\newcounter{ni}
\newcommand{\roii}{\roman{ni}}
\newcounter{san}
\newcommand{\roiii}{\roman{san}}
\newcounter{yon}
\newcommand{\roiv}{\roman{yon}}
\newcounter{go}
\newcounter{roku}
\newcounter{nana}
\newcommand{\Sc}{\mathscr{S}}
\newcommand{\Sp}{\mathscr{S}^\Bp}
\newcommand{\oSp}{\ol{\Sc}^\Bp}

\newcommand{\tSp}{\wt{\Sc}^{\Bp}}
\newcommand{\He}{\mathscr{H}}
\newcommand{\A}{\mathscr{A}}
\newcommand{\B}{\mathscr{B}}

\newtheorem{thm}{Theorem}[section]
\newtheorem{lem}[thm]{Lemma}
\newtheorem{cor}[thm]{Corollary}
\newtheorem{prop}[thm]{Proposition}

\def \para{\refstepcounter{thm} \par\medskip\noindent
                \textbf{\thethm .} }

\def \remark{\refstepcounter{thm} \par\medskip\noindent
                \textbf{Remark \thethm .} }

\def \remarks{\refstepcounter{thm} \par\medskip\noindent
                \textbf{Remarks \thethm .} }

\allowdisplaybreaks[4]
\numberwithin{equation}{thm}

\begin{document}
\setlength{\baselineskip}{4.9mm}
\setlength{\abovedisplayskip}{4.5mm}
\setlength{\belowdisplayskip}{4.5mm}
\renewcommand{\theenumi}{\roman{enumi}}
\renewcommand{\labelenumi}{(\theenumi)}
\renewcommand{\thefootnote}{\fnsymbol{footnote}}
\renewcommand{\thefootnote}{\fnsymbol{footnote}}
\parindent=20pt
\newcommand{\dis}{\displaystyle}

\medskip
\begin{center}
{\large \bf The representation type of 
Ariki-Koike algebras and 
cyclotomic $q$-Schur algebras} 
\\
\vspace{1cm}
Kentaro Wada\footnote{This research was  supported  by JSPS Research Fellowships for Young Scientists} 
\\ 
\vspace{0.5cm}
Graduate School of Mathematics \\
Nagoya University  \\
Chikusa-ku, Nagoya 464-8602,  Japan \\
(e-mail: kentaro-wada@math.nagoya-u.ac.jp) 

\end{center}
\title{}
\maketitle
\markboth{Kentaro Wada}{The representation type of Ariki-Koike algebras and cyclotomic $q$-Schur algebras}


\begin{abstract}
We give a necessary and sufficient condition 
on parameters  
for Ariki-Koike algebras (resp. cyclotomic $q$-Schur algebras)  
to be  of finite representation type. 
\end{abstract}


\setcounter{section}{-1}
\section{Introduction}

Let 
$F$ 
be an algebraically closed field, 
and 
$\CA$ 
be a finite dimensional associative algebra 
over 
$F$.  
We say that 
$\CA$  
is of finite representation type 
(simply, finite type) 
if there are only a finite number of isomorphism classes 
of 
indecomposable 
$\CA$-modules, 
and 
that 
$\CA$ 
is of infinite representation type (infinite type)  otherwise. 
Moreover, 
the infinite representation type has 
two classes, 
namely, 
tame type 
and 
wild type. 
$\CA$ is of tame type 
if 
indecomposable modules 
in each dimension 
come in one parameter families with finitely many exceptions. 
$\CA$ is of wild type 
if 
its module category is comparable with that of the free algebra in two variables. 
For precise definitions, see \cite{E-book} or \cite{B-book}. 
By Drozd's theorem, 
it is known that 
any finite dimensional algebra 
has 
finite type, 
tame type 
or 
wild type.

We consider the representation type 
of the Ariki-Koike algebra 
$\He_{n,r}=\He_{n,r}(q, \\ Q_1,\cdots,Q_r)$ 
over $F$  
with parameters 
$q,Q_1,\cdots,Q_r \in F$
and 
of 
the 
cyclotomic $q$-Schur algebra 
$\Sc_{n,r}=\Sc_{n,r}(q,Q_1,\cdots,Q_r)$ 
associated to 
$\He_{n,r}$.  
In the case where $r=1$, 
$\Sc_{n,1}$ 
is the $q$-Schur algebra, 
and 
the representation type 
of $\Sc_{n,1}$ 
has been determined by 
Erdmann-Nakano \cite{EN01}. 
On the other hand, 
the representation type 
of Hecke algebras of classical type 
has been determined by 
Uno \cite{Uno92}, 
Erdmann-Nakano \cite{EN01A}, 
Ariki-Mathas \cite{AM04} 
and 
Ariki \cite{Ari05}.  
In this paper, 
we will give 
a necessary and sufficient condition 
(here, we denote this condition by (CF))  
for 
$\He_{n,r}$ and 
$\Sc_{n,r}$ 
to be of finite type 
(Theorem \ref{main thm}). 

Suppose that $\Sc_{n,r}=\Sc_{n,r}(q,Q_1,\cdots,Q_r)$ satisfies the condition (CF). 
In order to show the finiteness,  
we will see that 
any block of 
$\Sc_{n,r}(q,Q_1,\cdots,Q_r)$ 
($r \geq 3$) 
is Morita equivalent to 
a certain block of 
$\Sc_{n',2}(q,Q_i,Q_{j})$ 
for some $i,j \in \{1,\cdots,r\}$ 
(Corollary \ref{Morita}).
Thus, 
the finiteness of 
$\Sc_{n,r}$ ($r\geq 3$)  
is reduced to the case where 
$r=2$. 
For the case where 
$r=2$, 
the finiteness is shown 
in a similar way as 
in the case of 
Hecke algebras of type B 
(\cite[Theorem 6.2]{AM04}). 
The finiteness for $\He_{n,r}$ 
follows from the finiteness for $\Sc_{n,r}$ 
(see Lemma \ref{lem Hecke Schur}). 

On the other hand, 
suppose that 
$\He_{n,r}$ 
does not satisfy the condition (CF). 
In order to show 
the infiniteness, 
we will make use of some properties of the structures of 
$\He_{n,r}$ which follows from the structures of 
$\Sc_{n,r}$ obtained 
in \cite{SW1} and \cite{W}. 
By adding some facts to results in \cite{SW1} and \cite{W}, 
we have the following picture (Theorem \ref{structure theorem}). 

\begin{picture}(30,75)
\put(0,55){$\dis \bigoplus_{\e=(n_1,\cdots,n_g)} \hspace{-1em} e_{\e}\Sc_{n,r}e_{\e}$ 
\scalebox{2}[1]{$\hookrightarrow $} 
$\Sp$ \scalebox{2}[1]{$\hookrightarrow $} $\tSp$ \scalebox{2}[1]{$\hookrightarrow $} $\Sc_{n,r}$}
\put(30,0){$ \dis \oSp \cong \hspace{-1em} \bigoplus_{\e=(n_1,\cdots,n_g)} \hspace{-1em} 
\Sc_{n_1,r_1}(q,\BQ_1) \otimes 
\Sc_{n_2,r_2}(q,\BQ_2) \otimes 
\cdots \otimes 
\Sc_{n_g,r_g}(q,\BQ_g) $}
\thicklines 
\put(35,12){\rotatebox{90}{\scalebox{2.2}[1.2]{$\twoheadleftarrow$}}}
\put(50,10){\rotatebox{30}{\scalebox{7}[0.6]{$\twoheadleftarrow$}}}
\put(60,10){\rotatebox{20}{\scalebox{10}[0.6]{$\twoheadleftarrow$}}}
\end{picture}\vspace{2em}\\
This implies 
the surjective homomorphism 
\[
\He_{n,r}(q,Q_1,\cdots,Q_r) \ra 
	\He_{n,k}(q,Q_{r-k+1},\cdots,Q_r)
\]
for $k=1,\cdots,r-1$ (Proposition \ref{structure theorem Hecke}).  
By using this surjection, 
the infiniteness of 
$\He_{n,r}$ 
is reduced to 
some special cases, 
namely 
the case where $r=2$ (in \cite{AM04}), 
Proposition \ref{infinite S3} 
and Proposition \ref{infinite condition g}. 
Finally, 
the infiniteness for $\Sc_{n,r}$ 
follows from the infiniteness for $\He_{n,r}$ 
(Lemma \ref{lem Hecke Schur}). 
\vspace{1em}\\
\textbf{Acknowledgment}\,\, 
The author is grateful to  
Professor S. Ariki for his suggestion on this problem. 
Thanks are also due to  
Professors T. Shoji and  H. Miyachi for many helpful advices and discussions. 
The author thanks Professor Shoji again 
for pointing out an error in a preliminary version of this paper.  
He also thanks the referee for many valuable comments. 

\section{The representation type of algebras }
Throughout this paper, 
we suppose that 
$F$ 
is an algebraically closed field, and any algebra 
$\CA$ 
is a finite dimensional unital associative algebra over 
$F$. 
We say just an 
$\CA$-module 
for a right 
$\CA$-module. 
The 
following results for the representation type are well known 
(see \cite{E-book}). 

\begin{lem}
\label{fundamental rep type}
Let 
$\CA$ 
be an 
$F$-algebra. 
\begin{enumerate}
\item 
Let 
$I$ 
be a two-sided ideal of 
$\CA$. 
If 
$\CA/I$ 
is of infinite (resp. wild) type then 
$\CA$ 
is also of  infinite (resp. wild) type.
\item 
Let 
$e$ 
be an idempotent of 
$\CA$. 
If 
$e\CA e$ 
is of infinite (resp. wild) type  then 
$\CA$ 
is also of infinite (resp. wild) type. 

\item
If an idempotent 
$e\in \CA$ 
is primitive then 
$\End_{\CA}(e\CA)\cong e \CA e$ 
is local. 
Moreover, 
$\End_{\CA}(e\CA)$ 
is of finite type 
if and only if 
$\End_{\CA}(e\CA)\cong F[x]/\lan x^m \ran$ 
for some integer 
$m\geq 0$, 
where 
$F[x]$ 
is a polynomial ring over 
$F$ 
with an indeterminate 
$x$, 
and 
$\lan x^m \ran $ 
is the ideal generated by the polynomial 
$x^m$.

\item
Let 
$P_1,\cdots, P_k$ 
be the complete set of non-isomorphic projective indecomposable 
$\CA$-modules. 
Then 
$\CA$ 
is Morita equivalent to 
$\End_{\CA}(P_1\oplus \cdots \oplus P_k)$. 
Thus if 
$\End_{\CA}(P_i)$ 
is of infinite type (resp. wild) for some 
$i$,  
then 
$\CA$ 
is of infinite (resp. wild) type. 
\end{enumerate}
\end{lem}


\para 
\textit{Cellular algebras.} 
A cyclotomic $q$-Schur algebra is a cellular algebra in the sense of 
\cite{GL96}. 
So, we give some fundamental properties of cellular algebras 
which we will be needed in later discussions. 
For more details for cellular algebras, see 
\cite{GL96} or 
\cite{M-book}. 

Let 
$\A$ 
be a cellular algebra over 
$F$ 
with respect to a poset 
$(\vL^+,\geq)$ 
and 
an algebra anti-automorphism 
$\ast$  
of 
$\A$.
Then we can define a cell module 
$W^\la$ 
for each 
$\la \in \vL^+$. 
Let 
$\rad W^\la$ 
be the radical of 
$W^\la$ 
with respect to the canonical bilinear form on 
$W^\la$. 
Put 
$L^\la =W^\la/\rad W^\la$. 
Since 
$\rad W^\la$ 
is an 
$\A$-submodule of 
$W^\la$, $L^\la$ 
is also an $\A$-module.
Set 
$\vL^+_0=\{\la \in \vL^+\,|\,L^{\la}\not=0\}$. 
Note that, for 
$\la\in \vL^+_0$, 
$\rad W^\la$ 
coincides with the Jacobson radical of 
$W^\la$. 
Then 
$\{L^\la\,|\,\la \in \vL^+_0\}$ 
is a complete set of non-isomorphic simple $\A$-modules. 
It is known that 
a cellular algebra $\CA$ 
is a quasi-hereditary algebra 
if and only if 
$\vL^+_0=\vL^+$. 

For 
$\la\in \vL^+$ 
and  
$\mu \in \vL^+_0$,  
let 
$d_{\la\mu}=[W^\la:L^\mu]$ 
be the decomposition number, 
namely 
the multiplicity of 
$L^\mu$ 
in 
the composition series of 
$W^\la$. 
If 
$d_{\la\mu}\not=0$ 
then 
$\la \geq \mu$.
The decomposition matrix of 
$\A$ 
is a matrix 
$D=(d_{\la\mu})_{\la\in \vL^+, \mu \in \vL^+_0}$. 
For 
$\la \in \vL^+_0$, 
let 
$P^\la$ 
be the projective cover of 
$L^\la$. 
The Cartan matrix of 
$\A$ 
is a matrix 
$C=(p_{\la\mu})_{\la,\mu \in \vL^+_0}$, 
where 
$p_{\la\mu}=\dim_F \Hom_{\A}(P^\la,P^\mu)=[P^{\mu}:L^\la]$. 
It is known that 
\begin{align}
\label{cartan}
C=\,^tDD.
\end{align}
Moreover, for 
$\la \in \vL^+_0$, 
$P^\la$ 
has a cell module filtration 
in which each cell module $W^\mu$ occurs with multiplicity $d_{\mu\la}$. 
The following properties  are well known 
(see \cite[2.5.]{EN01}). 

\begin{lem} 
\label{Ext}
Let 
$\A$ 
be a cellular algebra with respect to a poset 
$(\vL^+,\geq)$. 
\begin{enumerate}

\item 
$L^\la$ $(\la \in \vL^+_0)$ is self-dual. 
Thus we have,  
for 
$\la,\mu \in \vL^+_0$,   
\[\Ext_{\A}^i(L^\la,L^\mu)\cong \Ext^i_{\A}(L^\mu,L^\la) \text{ for any }i\geq 0.\]
In particular, 
\begin{align*}
\Ext_{\A}^1(L^\la,L^\mu) 
&\cong \Hom_{\A}(\rad P^{\la}/ \rad^2 P^{\la}, L^\mu)\\
&\cong \Hom_{\A}(\rad P^{\mu}/\rad^2 P^{\mu}, L^\la)\\
&\cong \Ext_{\A}^1(L^\mu,L^\la). 
\end{align*}

\item
If 
$\A$ 
is a quasi-hereditary algebra, 
then 
for 
$\la,\mu \in \vL^+$ 
such that 
$\la \geq \mu$, 
we have 
\begin{align*}
\Ext_{\A}^1(L^\la,L^\mu) &\cong \Hom_{\A}(\rad P^\la/\rad^2 P^\la, L^\mu)\\
	&\cong \Hom_{\A}(\rad W^\la/\rad^2 W^\la, L^\mu).
\end{align*}
Moreover,  
if
$\Ext_{\A}^1(L^\la, L^\mu) \not=0$ 
then 
$\la > \mu$ 
or 
$\mu > \la$.
\end{enumerate}

\end{lem} 

\para 
\textit{Tensor products and representation type.} 
For two cellular algebras  
$\A$ 
and  
$\B$, 
the tensor product 
$\A \otimes_F \B$
becomes a cellular algebra again 
in the natural way. 
For the representation type of 
$\A\ot_F \B$, 
we have the following lemma. 

\begin{lem}
\label{tensor of cellular}
Let $\A$, $\B$ be cellular algebras. 
\begin{enumerate}

\item
If $\B$ is semisimple, 
then 
the representation type of 
$\A\ot_F \B$ 
coincides with 
the representation type of 
$\A$. 

\item
If neither  $\A$ nor $\B$
are  semisimple, 
then 
$\A\ot_F\B$ is of infinite type. 

\item
If neither $\A$ nor $\B$
are  semisimple,  
and 
$\A$ contains a block 
which has at least 
three non-isomorphic simple modules 
as composition factors,  
then 
$\A\ot_F\B$ is of wild type.

\end{enumerate}
\end{lem}

\begin{proof}
(\roi) is clear. 
We show only (\roiii) 
since  
(\roii) is proven in a similar way. 

In order to apply 
\cite[Lemma 17]{Ari05} 
to 
$\A \otimes_F \B$, 
we consider the Gabriel quiver of 
$\A \otimes_F \B$.  
By assumption and Lemma \ref{Ext} (\roi), 
the Gabriel quiver of 
$\A$ 
contains the quiver 
\begin{picture}(70,0)
\put(0,0){$\bullet$}
\put(10,3){$\longrightarrow$}
\put(10,-3){$\longleftarrow$}
\put(30,0){$\bullet$}
\put(40,3){$\longrightarrow$}
\put(40,-3){$\longleftarrow$}
\put(60,0){$\bullet$}
\end{picture},
and 
the Gabriel quiver of 
$\B$ 
contains the quiver 
\begin{picture}(40,0)
\put(0,0){$\bullet$}
\put(10,3){$\longrightarrow$}
\put(10,-3){$\longleftarrow$}
\put(30,0){$\bullet$}
\end{picture}
as a subquiver. 
Thus, 
by 
\cite[Lemma 1.3]{Les94}, 
the Gabriel quiver of 
$\A\ot_F\B$ 
contains the quiver 

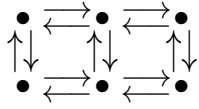
\hspace{5em}
\begin{picture}(40,45)
\put(3,26){$\bullet$}
\put(13,29){$\longrightarrow$}
\put(13,23){$\longleftarrow$}
\put(33,26){$\bullet$}
\put(43,29){$\longrightarrow$}
\put(43,23){$\longleftarrow$}
\put(63,26){$\bullet$}
\put(0,12){\scalebox{1}[1.5]{$\uparrow $}}
\put(6,12){\scalebox{1}[1.5]{$\downarrow $}}
\put(30,12){\scalebox{1}[1.5]{$\uparrow $}}
\put(36,12){\scalebox{1}[1.5]{$\downarrow $}}
\put(60,12){\scalebox{1}[1.5]{$\uparrow $}}
\put(66,12){\scalebox{1}[1.5]{$\downarrow $}}
\put(3,0){$\bullet$}
\put(13,3){$\longrightarrow$}
\put(13,-3){$\longleftarrow$}
\put(33,0){$\bullet$}
\put(43,3){$\longrightarrow$}
\put(43,-3){$\longleftarrow$}
\put(63,0){$\bullet$}
\end{picture}\vspace{3mm}\\
as subquiver. 
Thus 
$\A\ot_F\B$ 
is of wild type by 
\cite[Lemma 17]{Ari05}. 
\end{proof}

\para
\label{def Am}
Now we study a particular algebra defined by 
a quiver and relations. 
Let $m$ be a positive integer, 
and  
let 
$Q$ 
be a quiver 

\begin{picture}(100,30)
\put(10,10){$1 \qquad \,\,\,\,\,2 \qquad \,\,\,\,\,\,\cdots \qquad \quad \,\,\,m-1 \qquad\quad \, m,$}
\put(16,22){$ \underrightarrow{\,\,\,\,\,\a_{1}\,\,\,\,\,}$}
\put(16,0){$ \overleftarrow{\,\,\,\,\,\b_{1}\,\,\,\,\,}$}
\put(55,22){$ \underrightarrow{\,\,\,\,\,\a_{2}\,\,\,\,\,}$}
\put(55,0){$ \overleftarrow{\,\,\,\,\,\b_{2}\,\,\,\,\,}$}
\put(114,22){$ \underrightarrow{\,\,\a_{m-2}\,\,}$}
\put(114,0){$ \overleftarrow{\,\,\b_{m-2}\,\,}$}
\put(182,22){$ \underrightarrow{\,\,\a_{m-1}\,\,}$}
\put(182,0){$ \overleftarrow{\,\,\b_{m-1}\,\,}$}
\end{picture}\\
and 
$\CI$ 
be the two-sided ideal of the path algebra 
$FQ$ 
generated by the relations 
\begin{align*} 
\a_{m-1}\b_{m-1}=0, \,\,\,
\a_{i+1}\a_{i}=0, \,\,
\b_{i}\b_{i+1}=0,\,\,
\a_{i}\b_{i}=\b_{i+1}\a_{i+1} \,\,
\text{ for }1\leq i \leq m-2, 
\end{align*}
where we denote by 
$\a_{i+1}\a_{i}$ 
the path 
$i \stackrel{\a_{i}\,}{\longrightarrow } i+1 \stackrel{\a_{i+1}\,}{\longrightarrow } i+2$, 
etc.. 
We define 
$\CA_m=FQ/\CI$. 
Under the natural surjection 
$FQ \ra \CA_m$, 
we denote the image of paths 
in $FQ$ by the same symbol.  
By definition, 
$\CA_m$ 
has an 
$F$-basis 
\begin{align*}
e_1, \quad 
\begin{matrix} e_2, & \b_{1},\\ \a_{1}, & \a_{1}\b_{1}, \end{matrix} \quad\,\,
\begin{matrix} e_3, & \b_{2},\\ \a_{2}, & \a_{2}\b_{2}, \end{matrix} \quad
\cdots \quad
\begin{matrix} e_{m}, & \b_{m-1},\\ \a_{m-1}, & \a_{m-1}\b_{m-1}, \end{matrix}
\end{align*}
where $e_i$ is the path of length $0$ on the vertex $i$. 
It is known that 
$\CA_m$ 
is of finite type 
by \cite[3.1]{Erd91}. 
The following proposition 
was inspired by \cite[Proposition 3.2]{Erd91}, and 
will be used to prove 
the finiteness of cyclotomic $q$-Schur algebras. 

\begin{prop}
\label{decom matrix}
Let 
$\A$ 
be a cellular algebra with respect to the poset 
$(\vL^+,\geq)$. 
If 
$\A$ 
is a quasi-hereditary algebra 
with 
decomposition matrix 
\[D=\left(
\begin{matrix}
1&&&\\
1&1&&\\
&\hspace{-1mm}\ddots\hspace{-1mm}&\hspace{-1mm}\ddots\hspace{-1mm}&\\
&&1&1
\end{matrix}
\right)
\quad
\text{(all omitted entries are zero)},
\]
and 
if any projective indecomposable $\A$-module has the simple socle, 
then
$\A$ 
is Morita equivalent to 
$\CA_m$ 
with 
$m=|\vL^+|$. 
In particular, 
$\A$ 
is of finite type. 
\end{prop}
\begin{proof}

Let 
$\vL^+=\{\la_1,\cdots,\la_m\}$ 
such that 
$i<j$ 
if 
$\la_i<\la_j$. 
We denote the simple $\A$-module 
by 
$L^{\la_i}$, 
and its projective cover 
by 
$P^{\la_i}$ 
for 
$\la_i \in \vL^+$.
By (\ref{cartan}), 
the Cartan matrix of 
$\A$ 
is 
\[C=\,^tDD=
\left(
\begin{matrix}
2&1&&&&&\\
1&2&1&&&&\\
&1&2&1&&&\\
&&\hspace{-1mm}\ddots\hspace{-1mm}&\hspace{-1mm}\ddots\hspace{-1mm}&\hspace{-1mm}\ddots\hspace{-1mm}&&\\
&&&1&2&1\\
&&&&1&1
\end{matrix}
\right)
\quad 
\text{(all omitted entries are zero)}.
\]
Combined with the second isomorphism in  Lemma \ref{Ext} (\roii), 
we have 
\begin{align}
\label{radical layer}
P^{\la_1}=
\begin{matrix}
L^{\la_1}\\
L^{\la_2}\\
L^{\la_1} 
\end{matrix} 
\,\,, \qquad
P^{\la_i}=
\begin{matrix}
L^{\la_i}\\
L^{\la_{i-1}} \oplus L^{\la_{i+1}}\\
L^{\la_i}
\end{matrix}
\quad \text{for }
i=2,\cdots,m-1\,, \quad 
P^{\la_m}=
\begin{matrix}
L^{\la_m}\\
L^{\la_{m-1}}
\end{matrix}, 
\end{align}
where the $i$th row in the right hand side of equation 
means the $i$th radical layer of 
$P^{\la_i}$. 
We also have 
\begin{align}
\label{dimension}
\dim_F \Hom_{\A}(P^{\la_i},P^{\la_j})
=
\begin{cases}
1 &\text{ if } |i-j|=1 \text{ or } i=j=m, \\
2 &\text{ if }i=j=1,\cdots,m-1,\\
0 &\text{ if }|i-j|>1.
\end{cases}
\end{align}
Let 
$\A'=\End_{\A}\big(\bigoplus_{i=1}^{m}P^{\la_i}\big)$. 
Then 
$\A$ 
is Morita equivalent to 
$\A'$. 
Let 
$e_i' \in \A'$ 
be the identity map on 
$P^{\la_i}$ 
and be the 
$0$-map on 
$P^{\la_j}$ 
($i\not=j$) 
for 
$i=1,\cdots ,m$. 

For $i=2,\cdots,m$, 
let $M^{\la_i}$ and $N^{\la_i}$    
be the $\A$-modules such that 
\[ 
M^{\la_i}=\begin{matrix} L^{\la_{i-1}} \\ L^{\la_i} \end{matrix}, 
\qquad 
N^{\la_i}=\begin{matrix} L^{\la_{i}} \\ L^{\la_{i-1}} \end{matrix}.
\] 
Note that 
any projective indecomposable $\CA$-module has the simple socle.  
Then,   
by \eqref{radical layer}, 
one can take the natural surjective homomorphism 
$\vf_i : P^{\la_i} \ra M^{\la_{i+1}}$ for $i=1,\cdots, m-2$,  
and take the injective homomorphism 
$\psi_i: M^{\la_{i}} \ra P^{\la_i}$ for $i=2,\cdots, m-1$. 
Put 
$\a'_i = \psi_{i+1} \circ \vf_i \in \Hom_{\CA}(P^{\la_i}, P^{\la_{i+1}})$ for $i=1,\cdots,m-2$. 
We also define 
$\a'_{m-1} \in \Hom_{\CA}(P^{\la_{m-1}}, P^{\la_m})$ 
by the composition of 
the natural surjection $P^{\la_{m-1}} \ra L^{\la_{m-1}}$ 
and the injection $L^{\la_{m-1}} \ra P^{\la_m}$. 
Similarly, 
one can define 
$\b_i' \in \Hom_{\CA}(P^{\la_{i+1}},P^{\la_i})$ for $i=1,\cdots,m-1$ 
such that 
$\Im \b_i' = N^{\la_{i+1}}$. 
We regard 
$\a_i'$ 
(resp. $\b_i'$) 
as an element of 
$\A'$ 
by 
$\a_i'(P^{\la_j})=0$ 
(resp. $\b_i'(P^{\la_j})=0$) 
for any 
$j\not=i$ 
(resp. $j\not=i+1$). 
Then 
we have 
$\Im \b_{i}'\a_i'=L^{\la_{i}}$ 
for $i=1,\cdots,m-1$, 
$\Im \a_{i}'\b_{i}'=L^{\la_{i+1}}$ 
for $i=1,\cdots,m-2$ 
and  
$\Im \a_{m-1}'\b_{m-1}'=0$. 
Since 
$\dim_F \Hom_{\A}(P^{\la_i}, L^{\la_i})=1$, 
we have 
$\a_{i}'\b_{i}'=\b_{i+1}'\a_{i+1}'$ 
for 
$i=1,\cdots,m-2$ 
by multiplying $\a_i'$  by a scalar  
if  necessary. 
Moreover, we have 
$\a_{i+1}'\a_{i}'=0$ 
and 
$\b_{i}'\b_{i+1}'=0$ 
since 
$\Hom_{\A}(P^{\la_i},P^{\la_j})=0$ 
for 
$|i-j|>1$. 
Now we can define a surjective  homomorphism of algebras 
from 
$\CA_m$ 
to 
$\A'$ 
by 
$X \mapsto X'$ 
($X\in \{ e_i,\, e_m,\, \a_i,\,\b_i\,|\, i=1,\cdots,m-1\}$), 
and 
we see that 
this gives an isomorphism 
by comparing the dimensions. 
\end{proof}

The following lemma will be used to prove 
the infiniteness of Ariki-Koike algebras.

\begin{lem}
\label{lem infinite cellular}
Let $\CA$ be a cellular algebra with respect to the 
poset $\vL^+=\{\la_0,\la_1,\la_2\}$. 
If $\CA$ has the following decomposition matrix (thus, $\CA$ is a quasi-hereditary) 
\\

\begin{tabular}{c|ccc}

&$L^{\la_0}$ &$L^{\la_1}$ &$L^{\la_2}$ \\ \hline 
$W^{\la_0}$&1&0&0\\
$W^{\la_1}$&a&1&0\\
$W^{\la_2}$&b&0&1 
\end{tabular}
\quad $( a,b >0)$, 
\\[2mm]
then 
$\CA$ is of infinite type. 
\end{lem}

\begin{proof} 
From the decomposition matrix, 
we see that 
$W^{\la_1}$ (resp. $W^{\la_2}$) 
has the unique simple top 
$L^{\la_1}$ (resp. $L^{\la_2}$),   
and that 
any other composition factor  
of $W^{\la_1}$ (resp. $W^{\la_2}$)  
is isomorphic to $L^{\la_0}$. 
By the general theory of cellular algebras, 
$P^{\la_0}$ has the filtration 
$ P^{\la_0} =P_0 \supsetneq P_1 \supsetneq P_2 \supsetneq 0$ 
such that 
$P_0/P_1 \cong W^{\la_0} \cong L^{\la_0}$, 
$P_1/ P_2 \cong (W^{\la_1})^{\oplus a}$ 
and  
$P_2 \cong (W^{\la_2})^{\oplus b}$. 
This filtration implies that  
$L^{\la_1}$ appears exactly $a$ times in the second radical layer of $P^{\la_0}$. 
On the other hand, 
by Lemma \ref{Ext}, 
we see that 
$L^{\la_2}$ appears at least once in the  second radical layer of $P^{\la_0}$. 
Thus, 
$L^{\la_0}$ appears at least twice in the third radical layer of $P^{\la_0}$. 
This implies that 
$\End_{\CA}(P^{\la_0}/ \rad^4 P^{\la_0})$ 
is not isomorphic to 
$F[x]/ \lan x^m \ran$ for any $m \geq 0$. 
Combining  with Lemma \ref{fundamental rep type}, 
we have that 
$\CA / \rad^4 \CA$ 
is of infinite type, 
thus 
$\CA$ is also of infinite type.

\end{proof}

\section{Ariki-Koike algebras and Cyclotomic $q$-Schur algebras}
In this section, 
we introduce  Ariki-Koike algebras and  cyclotomic $q$-Schur algebras.   
We also give   
some properties of them.  
\para
\label{comb}
A composition
$\mu=(\mu_1,\mu_2,\cdots)$ 
is a finite sequence of non-negative integers, 
and 
$|\mu|=\sum_{i} \mu_i$ 
is called the size of 
$\mu$. 
If 
$\mu_l\not=0$ 
and 
$\mu_k=0$ 
for any 
$k>l$, 
then 
$l=\ell(\mu)$ is called the length of 
$\mu$. 
If the composition 
$\la$ 
is a weakly decreasing sequence, 
$\la$ 
is called a partition. 
An 
$r$-tuple 
$\mu=(\mu^{(1)},\cdots ,\mu^{(r)})$ 
of compositions is called an $r$-composition, 
and the size 
$|\mu|$ 
of 
$\mu$ 
is  defined by 
$|\mu|=\sum_{i=1}^r|\mu^{(i)}|$. 
In particular, 
if all 
$\mu^{(i)}$ 
are partitions, 
$\mu$ 
is called an 
$r$-partition. 
For 
$n,r\in \ZZ_{>0}$ 
and 
an $r$-tuple 
$\Bm=(m_1,\cdots ,m_r)\in \ZZ_{>0}^r$, 
we denote by 
$\vL_{n,r}(\Bm)$ 
the set of 
$r$-compositions 
$\mu=(\mu^{(1)},\cdots ,\mu^{(r)})$ 
such that 
$|\mu|=n$ 
and 
that 
$\ell(\mu^{(k)}) \leq m_k$ 
for 
$k=1,\cdots ,r$. 
We define 
$\vL^+_{n,r}(\Bm)$ 
as the subset of 
$\vL$ 
consisting of 
$r$-partitions. 
Throughout this paper, 
we assume the following condition for $\vL_{n,r}(\Bm)$. 
\begin{description}
\item[(CL)] $m_i \geq n$ for any $i=1,\cdots, r$. 
\end{description}
Under this condition, 
$\vL^+_{n,r}(\Bm)$ 
coincides with  
the set of 
$r$-partitions 
of size $n$. 
In particular, 
$\vL^+_{n,r}(\Bm)$ 
is independent of  
a choice of 
$\Bm$ 
satisfying (CL). 
Thus, 
we write it 
simply as   
$\vL_{n,r}^+$ 
instead of 
$\vL_{n,r}^+(\Bm)$.  
Similarly, 
we may write 
$\vL_{n,r}(\Bm)$ 
simply as 
$\vL_{n,r}$ 
if 
there is no fear of confusion. 

For 
$\mu \in \vL_{n,r}(\Bm)$, 
the diagram of 
$\mu$ 
is the set  
\[ 
[\mu] = \big\{(i,j,k) \in  \NN \times \NN \times \{1,2,\cdots,r \} 
\bigm| 1 \leq j \leq \mu_i^{(k)}\big\}.
\]
We call an element of $[\mu]$ a node. 

We define a partial order, 
the so-called 
\lq\lq dominance order", 
on 
$\vL_{n,r}(\Bm)$ 
by 
$\mu \trreq \nu$ 
if and only if
\[\sum_{i=1}^{l-1}|\mu^{(i)}| +\sum_{j=1}^{k}\mu^{(l)}_j \geq \sum_{i=1}^{l-1}|\nu^{(i)}| +\sum_{j=1}^{k}\nu^{(l)}_j\]
for any 
$1\leq l \leq r$, $1\leq k \leq m_l$. 
If 
$\mu \trreq \nu $ 
and 
$\mu \not= \nu$, 
we write it as 
$\mu \trr \nu$.

\para
\textit{Ariki-Koike algebras.} 
\label{definition of Ariki-Koike}
Let 
$F$ 
be an algebraically closed field, 
and take 
$q\not=0,Q_1,\cdots,Q_r \in F$.
The Ariki-Koike algebra 
$\He_{n,r}=\He_{n,r}(q,Q_1,\cdots,Q_r)$ 
is an associative algebra 
over 
$F$ 
with generators 
$T_0,T_1,\cdots,T_{n-1}$
and  
relations 
\begin{align*}
&(T_0-Q_1)(T_0-Q_2)\cdots (T_0-Q_r)=0,\\
&T_0 T_1 T_0 T_1=T_1 T_0 T_1 T_0,\\
&(T_i+1)(T_i-q)=0 \quad \text{ for }i=1,\cdots, n-1,\\ 
&T_i T_{i+1}T_i=T_{i+1}T_i T_{i+1} \quad \text{ for }i=1,\cdots, n-2,\\
&T_i T_j=T_j T_i \quad \text{ if }|i-j|\geq 2.   
\end{align*}

By \cite{DJM98}, 
it is known that 
$\He_{n,r}$ is a cellular algebra with a cellular basis 
\[ 
\big\{ m_{\Fs \Ft} \bigm| \Fs,\Ft \in \Std(\la) \text{ for some } \la \in \vL_{n,r}^+ \big\},
\]
where 
$\Std(\la)$ is the set of standard tableaux of shape $\la$ 
(see \cite{DJM98} for definitions).  

We denote by 
$S^\la$ 
the cell module of $\He_{n,r}$ 
corresponding to $\la \in \vL_{n,r}^+$, 
which is called the Specht module. 
Put 
$D^\la = S^\la / \rad S^\la$, 
where 
$\rad S^\la$ is the radical of $S^\la$ 
with respect to the canonical bilinear form on $S^\la$. 
When 
$q \not=1$ and $ Q_i \not=0$ for $1 \leq i \leq r$, 
it is known that 
$D^\la \not= 0 $  
if and only if 
$\la $ is a Kleshchev multipartition 
by \cite{Ari01}.  
Thus 
$\{ D^\la \,|\, \la \in \vL_{n,r}^+ : \text{ Kleshchev multipartition} \}$ 
gives a complete set of non-isomorphic simple $\He_{n,r}$-modules 
(see \cite{Ari01} or \cite[\S 3.4]{Mat04}  for the definition of Kleshchev multipartitions and more details).

\para
\textit{Cyclotomic $q$-Schur algebras.}  
The cyclotomic $q$-Schur algebra 
$\Sc_{n,r}= \Sc_{n,r}(q,Q_1,\\ \cdots,Q_r)$ 
associated to 
$\He_{n,r} (q,Q_1,\cdots,Q_r)$ 
is defined by 
\[
\Sc_{n,r}
=\Sc_{n,r}(\vL_{n,r}(\Bm)) 
=\End_{\He_{n,r}}\Big(\bigoplus_{\mu \in \vL_{n,r}(\Bm )} M^\mu \Big), 
\]
where 
$M^\mu$ 
is a certain 
$\He_{n,r}$-module 
introduced in \cite{DJM98}  
with respect to 
$\mu \in \vL_{n,r}(\Bm)$.  
By [loc. cit.], 
$\Sc_{n,r}$ 
is a cellular algebra with respect to the poset 
$(\vL^+_{n,r},\trreq)$. 
More precisely, 
$\Sc_{n,r}$ 
has a cellular basis 
\[\ZC=\big\{ \vf_{ST} \bigm| S,T \in \CT_0(\la) \text{ for some } \la \in \vL^+_{n,r} \big\},\]
where 
$\CT_0(\la)=\bigcup_{\mu \in \vL_{n,r}}\CT_0(\la,\mu)$, 
and 
$\CT_0(\la,\mu)$ 
is the set of semistandard tableaux of shape 
$\la$ 
with type 
$\mu$
(see \cite{DJM98} for definitions). 
By definition, 
for 
$S\in \CT_0(\la,\mu)$ 
and 
$T \in \CT_0(\la,\nu)$, 
we have that 
$\vf_{ST} \in \Hom_{\He_{n,r}}(M^\nu,M^\mu)$ 
and 
$\vf_{ST}|_{M^\kappa}=0$ 
($\kappa \not= \nu$). 
For 
$\mu \in \vL_{n,r}$, 
let 
$\vf_{\mu} \in \Sc_{n,r}$ 
be the identity map on 
$M^\mu$ 
and zero-map on 
$M^{\kappa}$ ($\kappa \not=\mu$). 
Then 
we have 
\[1_{\Sc_{n,r}}=\sum_{\mu \in \vL_{n,r}} \vf_\mu,\]
and 
$\{\vf_\mu\,|\, \mu\in \vL_{n,r}\}$  
is a set of pairwise orthogonal idempotents. 
Thus, 
for 
$S\in \CT_0(\la,\mu),\,T\in \CT_0(\la,\nu)$ 
and 
$\kappa \in \vL_{n,r}$, 
 we have 
\begin{align}
\label{cellular basis of Schur}
\vf_{\kappa}\vf_{ST}=\d_{\kappa \mu}\vf_{ST}, \quad \vf_{ST}\vf_{\kappa}=\d_{\kappa \nu}\vf_{ST}, 
\end{align} 
where 
$\d_{\kappa \mu}=1$ 
if $\kappa = \mu$,  
and 
$\d_{\kappa \mu}=0$ 
if $\kappa \not= \mu$. 

Let 
$W^\la$ 
be the cell module 
for 
$\Sc_{n,r}$ 
corresponding to 
$\la \in \vL_{n,r}^+$,  
which is called 
the Weyl module, 
and 
$\rad W^\la$ 
be the radical of 
$W^\la$ 
with respect to the canonical bilinear form on 
$W^\la$. 
Put 
$L^\la=W^\la/\rad W^\la$. 
It is known that 
$L^\la \not=0$ 
for any 
$\la \in \vL_{n,r}^+$, 
namely 
the cyclotomic $q$-Schur algebra 
$\Sc_{n,r}$ 
is a quasi-hereditary algebra. 
Thus 
$\{L^\la\,|\,\la\in \vL_{n,r}^+\}$ 
is a complete set of non-isomorphic 
simple $\Sc_{n,r}$-modules. 

\begin{remarks}
\label{remark}\\
(\roi) 
By a general theory, 
for a cellular algebra $\A$ over any field $F$, 
$F$ is a splitting field for $\A$. 
Thus, 
we may assume  that 
$F$ is an algebraically closed field 
without loss of generality. \vspace{2mm}\\
(\roii)
A cyclotomic $q$-Schur algebra 
$\Sc_{n,r}$ 
depends on a choice of $\vL_{n,r}(\Bm)$. 
But, it is known that 
$\Sc_{n,r}(\vL_{n,r}(\Bm))$ 
is Morita equivalent to 
$\Sc_{n,r}(\vL_{n,r}(\Bm'))$ 
with same parameters 
if both of 
$\Bm$ 
and 
$\Bm'$ 
satisfy the condition (CL). \vspace{2mm}\\
(\roiii)
By definitions, 
we have the following properties. 
\begin{description}
\item[(a)]
For any $0 \not= c \in F$, 
we have an isomorphism  
$\He_{n,r}(q,Q_1,\cdots,Q_r) 
\cong 
\He_{n,r}(q, cQ_{1},\cdots, cQ_{r})$.  
We also have an isomorphism  
$\Sc_{n,r}(q, Q_1,\cdots,Q_r) \cong \Sc_{n,r}(q,cQ_{1},\cdots, cQ_{r})$. 

\item[(b)] 
For any permutation 
$\s$ 
of $r$ letters, 
we have an isomorphism  
$\He_{n,r}(q,Q_1,\\ \cdots, Q_r) 
\cong 
\He_{n,r}(q,Q_{\s(1)},\cdots, Q_{\s(r)})$.  
However, 
$\Sc_{n,r}(q, Q_1,\cdots,Q_r)$ 
is not isomorphic to $\Sc_{n,r}(q,Q_{\s(1)},\cdots, Q_{\s(r)})$ 
in general. 

\end{description}
\end{remarks}\quad\vspace{-1.5em}\\

\para  
Put $\omega=(-,\cdots,-,(1^n))$, 
where \lq\lq $-$" means the empty partition, 
then $\omega \in \vL_{n,r}^+$ 
by condition (CL) (see \ref{comb}).  
By the definition, 
$\vf_\omega$ 
is the identity map on 
$M^\omega$ 
and 
$0$-map on 
$M^{\kappa }$ ($\kappa \not=\omega$). 
In particular, 
$\vf_\omega$ 
is an idempotent of 
$\Sc_{n,r}$. 
It is well known that 
the subalgebra 
$\vf_\omega \Sc_{n,r} \vf_{\omega}$ 
of 
$\Sc_{n,r}(q,Q_1,\cdots,Q_r)$ 
is isomorphic to 
$\He_{n,r}(q,Q_1,\cdots,Q_r)$ 
as algebras. 
Thus, 
by Lemma \ref{fundamental rep type} (\roii), 
we have the following lemma. 

\begin{lem}

\label{lem Hecke Schur}
If 
$\He_{n,r}(q,Q_1,\cdots,Q_r)$ 
is of infinite (resp. wild) type 
then 
$\Sc_{n,r}(q,Q_1,\\\cdots,Q_r)$ 
is of infinite (resp. wild) type.  
\end{lem}

\para 
\label{Schur functor}
\textit{Schur functor.}  
Since $\He_{n,r}$ 
is isomorphic to the subalgebra 
$\vf_\omega \Sc_{n,r} \vf_\omega$ 
of  
$\Sc_{n,r}$, 
we can define a functor  
$F=\Hom_{\Sc_{n,r}}(\vf_\omega \Sc_{n,r}, -)$ 
from the category of finite dimensional $\Sc_{n,r}$-modules 
to the category of finite dimensional $\He_{n,r}$-modules. 
The following lemma is known 
(see e.g. \cite[Theorem 5.1]{Mat04}, \cite[Appendix]{Don}).

\begin{lem}\
\label{lem Schur functor}

\begin{enumerate}
\item 
$F(W^\la) \cong S^\la$ as $\He_{n,r}$-modules for $\la \in \vL^+_{n,r}$. 
\item 
$F(L^\la) \cong D^\la$ as $\He_{n,r}$-modules for $\la \in \vL^+_{n,r}$.
\item 
For $\la \in \vL_{n,r}^+$,  
let 
$P^\la$ 
be the projective cover of 
$L^\la$. 
Then we have that 
$\{ F(P^\la)\,|\, \la \in \vL^+_{n,r} \text{ such that } F(L^\la) \not=0 \}$ 
gives a complete set of non-isomorphic projective indecomposable $\He_{n,r}$-modules. 
\end{enumerate}
\end{lem}

\para
For later discussions, 
we describe some structural properties of 
$\Sc_{n,r}$, 
the essential part of which 
has been obtained in 
\cite{SW1} 
and 
\cite{W}. 
Here, 
we modify such results 
for our purpose. 
Take and fix 
$\Bp=(r_1,\cdots,r_g) \in \ZZ_{>0}^g$ 
such that 
$r_1+\cdots +r_g=r$. 
Put 
$p_i=\sum_{j=1}^{i-1} r_j$ 
with 
$p_1=0$. 
For 
$\mu \in \vL_{n,r}$, 
put 
$\mu^{[k]}=(\mu^{(p_k+1)},\mu^{(p_k+2)},\cdots,\mu^{(p_k+r_k)})$ 
for 
$k=1,\cdots, g$.
Thus, 
$\mu^{[k]}$ 
is a $r_k$-composition. 
We define a map 
$\a_{\Bp} : \vL_{n,r} \ra \ZZ_{\geq 0}^g$ 
by  
$\mu \mapsto (|\mu^{[1]}|,|\mu^{[2]}|,\cdots,|\mu^{[g]}|)$. 
Thus, 
we have 
\[
\Im \a_{\Bp}=\Delta_{n,g}:=\{(n_1,\cdots,n_g) \in \ZZ_{\geq 0}^g \,|\, n_1+\cdots+n_g=n\}.
\] 

Recall that, 
for 
$\la \in \vL^+$,  
$\CT_0(\la)=\bigcup_{\mu \in \vL_{n,r}}\CT_0(\la,\mu)$ 
is the set of semistandard tableaux of shape 
$\la$.  
We define two subsets 
$\CT_\Bp^+(\la)$ 
and 
$\CT_{\Bp}^-(\la)$ 
of 
$\CT_0(\la)$ 
by 
\[ 
\CT_{\Bp}^+(\la)=\bigcup_{\mu \in \vL_{n,r} \atop \a_{\Bp}(\la)=\a_{\Bp}(\mu)}\CT_0(\la,\mu), 
\quad 
\CT_\Bp^-=\CT_0(\la) \setminus \CT_\Bp^+(\la).
\]
Moreover, 
for each 
$\e \in \D_{n,g}$, 
we define the subset 
$\CT_\Bp^\e(\la)$ 
of 
$\CT_0(\la)$ 
by 
\[\CT_\Bp^\e (\la)=\bigcup_{\mu \in \vL_{n,r} \atop \a_{\Bp}(\mu)=\e}\CT_0(\la,\mu).\]
By definition, 
we have 

\begin{align}
\label{CTe}
\begin{cases}
\CT_\Bp^\e(\la)
= \CT_\Bp^+(\la) \quad \text{if }\a_{\Bp}(\la)=\e, \\ 
\CT_\Bp^\e(\la) \subseteq \CT_\Bp^-(\la) \quad \text{if }\a_{\Bp}(\la) \not=\e.
\end{cases}
\end{align}

For 
$\e \in \D_{n,g}$, 
put 
\[e_\e=\sum_{\mu \in \vL_{n,r} \atop \a_\Bp(\mu)=\e} \vf_\mu.
\]
Since 
$\{ \vf_\mu \,|\, \mu \in \vL_{n,r}\}$ 
is a set of pairwise orthogonal idempotents, 
$e_\e$ 
is also an idempotent. 
Thus, 
$\Sc^\e=e_\e \Sc_{n,r} e_\e$ 
is a subalgebra of 
$\Sc_{n,r}$. 
The following theorem 
is a modification of the results in \cite{SW1}, \cite{W}. 

\begin{thm}
\label{structure theorem}
For each 
$\e \in \D_{n,g}$, 
we have the following. 

\begin{enumerate}
\item 
$\Sc^\e$ 
is a subalgebra of 
$\Sc_{n,r}$. 
Moreover, 
$\Sc^\e$ 
is a cellular algebra 
with 
a cellular basis  
\[
\ZC^\e =\big\{\vf_{ST} \bigm| S,T \in \CT_\Bp^\e(\la) \text{ for some }\la \in \vL^+_{n,r}\}.
\]

\item 
Let  
$\wh{\Sc}^\e$ 
be the $F$-subspace of 
$\Sc^\e$ 
spanned by 
$\{\vf_{ST} \,|\, S,T\in \CT_\Bp^\e(\la) \text{ for }\la \in \vL_{n,r}^+ \text{ such that }\a_{\Bp}(\la)\not=\e\}$. 
Then 
$\wh{\Sc}^\e$ 
is a two sided ideal of 
$\Sc^\e$. 
Thus one can define a quotient algebra 
$\ol{\Sc}^\e=\Sc^\e/ \wh{\Sc}^\e$. 

\item
$\ol{\Sc}^\e$ 
is a cellular algebra with a cellular basis 
\[
\ol{\ZC}^\e= 
\{\ol{\vf}_{ST}\,|\, S,T \in \CT_\Bp^+(\la) \text{ for }\la\in \vL^+_{n,r} \text{ such that }\a_{\Bp}(\la)=\e\}.
\] 

\item 
There exists an isomorphism of algebras 
\[
\ol{\Sc}^\e \cong 
\Sc_{n_1,r_1}(q,\BQ_1) \otimes 
\Sc_{n_2,r_2}(q,\BQ_2) \otimes 
\cdots \otimes 
\Sc_{n_g,r_g}(q,\BQ_g), 
\]
where 
$\e=(n_1,\cdots,n_g)$ 
and 
$\BQ_k=\{Q_{p_k+1},\cdots,Q_{p_k+r_k}\}$ 
for $k=1,\cdots,g$. 

\end{enumerate}
\end{thm}

\begin{proof}
It is already shown that 
$\Sc^\e$ 
is a subalgebra of 
$\Sc_{n,r}$. 
By (\ref{cellular basis of Schur}), 
we see that 
$\ZC^\e$ 
is a basis of 
$\Sc^\e$. 
Thus, 
$\Sc^\e$ 
inherits 
the cellular structure 
from 
$\Sc_{n,r}$. 
This proves (\roi). 
By noting (\ref{CTe}),  
(\roii) follows from 
\cite[Lemma 2.11]{W}.  
Now  
(\roiii) is clear. 
(\roiv) 
follows from the proof of 
\cite[Theorem 4.15]{SW1}. 
\end{proof}

\remark 
The statements in 
Theorem \ref{structure theorem} 
except (\roiv)  
also hold  
for certain types cellular algebras 
under the setting in \cite{W}. 
\vspace{-0.5em}\\

Theorem \ref{structure theorem} 
implies the following proposition 
for Ariki-Koike algebras.

\begin{prop}
\label{structure theorem Hecke}
For $k=1,2,\cdots,r-1$,  
there exists a surjective homomorphism 
\[ 
\He_{n,r}(q,Q_1,\cdots,Q_r) 
\ra 
\He_{n,k} (q,Q_{r-k+1},\cdots, Q_{r-1}, Q_{r}). 
\]
\end{prop}

\begin{proof}
Put  
$\Bp=(r-k,k)$ 
and 
$\eta=(0,n)$.  
Then, by Theorem \ref{structure theorem}, 
we have a surjective homomorphism 
$\Sc^{\eta} \ra \ol{\Sc}^\eta$, 
where $\ol{\Sc}^\eta$ 
is isomorphic to 
$1 \otimes \Sc_{n,k} (q,Q_{r-k+1},\cdots, Q_r)$. 
Since 
$\a_{\Bp}(\omega) = \eta$, 
we have 
$\vf_{\omega} e_{\eta} = e_{\eta} \vf_{\omega} = \vf_\omega$. 
This implies that 
$\vf_{\omega} \Sc^\eta \vf_{\omega} \cong \He_{n,r}(q,Q_1,\cdots,Q_r)$. 
On the other hand, 
one see easily that 
\[
\ol{\vf_\omega} \ol{\Sc}^\eta \ol{\vf_\omega} \cong \He_{n,k}(q,Q_{r-k+1},\cdots,Q_r) 
\]
through the isomorphism  
$\ol{\Sc}^\eta \cong 1 \otimes \Sc_{n,k} (q,Q_{r-k+1},\cdots, Q_r)$, 
where $\ol{\vf_\omega}$ is the image of $\vf_\omega$ 
under the surjection $\Sc^\eta \ra \ol{\Sc}^\eta$. 
Hence 
the surjection 
$\Sc^\eta \ra \ol{\Sc}^\eta$ 
implies the proposition. 
\end{proof}

The following corollary plays an important role in later discussions. 

\begin{cor}
\label{cor induction}
Take 
a subset 
$\{i_1,\cdots,i_k\} \subset \{1, \cdots,r\}$. 
If 
$\He_{n,k}(q,Q_{i_1},\cdots,Q_{i_k})$
is of infinite (resp. wild) type, 
then 
$\He_{n,r}(q,Q_1,\cdots,Q_r)$ 
is also of infinite (resp. wild) type.  
\end{cor}

\begin{proof}
By Remarks \ref{remark} (\roiii), 
we may suppose that $\{i_1,i_2,\cdots,i_k\}=\{r-k+1,\cdots,r-1,r\}$. 
Hence,  
the corollary follows from 
Proposition \ref{structure theorem Hecke} 
together with Lemma \ref{fundamental rep type} (\roi). 
\end{proof}


\section{The representation type of Ariki-Koike algebras and cyclotomic $q$-Schur algebras}
In this section, 
we study the representation type of Ariki-Koike algebras and cyclotomic $q$-Schur algebras. 
First, we recall a necessary and sufficient condition 
for Ariki-Koike algebras (resp. cyclotomic $q$-Schur algebras)  
to be semisimple. 
By Ariki \cite{Ari94}, 
the condition for Ariki-Koike algebras  
to be semisimple 
was obtained, 
and through the double centralizer property (\cite[Theorem 5.3]{Mat04}), 
we have the following theorem.

\begin{thm}
$\He_{n,r}(q,Q_1,\cdots, Q_r)$ 
(resp. $\Sc_{n,r}(q,Q_1,\cdots,Q_r)$ ) 
is semisimple 
if and only if 
\[\prod_{i=1}^n (1+q+\cdots +q^{i-1}) \prod_{1 \leq i<j \leq r} \prod_{-n<a<-n}(q^a Q_i-Q_j) \not=0.\]
\end{thm}


\para 
In order to compute the decomposition numbers 
of $\He_{n,r}$ or $\Sc_{n,r}$ 
in some cases, 
we will use the Jantzen sum formula obtained by James-Mathas 
\cite{JM00}. 
Here, we review on the Jantzen sum formula briefly 
(see [loc. cit.] for more details). 

Let $R$ be a discrete valuation ring with the unique maximal ideal $\wp$,     
$K$ be the quotient field of $R$,   
and 
$F$ be the residue field $R/ \wp$.  
Let 
$\nu_{\wp} : R^\times \ra \NN$ 
be the 
$\wp$-adic valuation map, 
and 
we extend 
$\nu_{\wp}$ 
to a map 
$K^\times \ra \ZZ$ 
in the natural way.   

Take 
$\wh{q}, \wh{Q}_1,\cdots ,\wh{Q}_r \in R$, 
where $\wh{q}$ is invertible in $R$.
Let 
$q, Q_1,\cdots Q_r \in F$ 
be the image of 
$\wh{q}, \wh{Q}_1,\cdots ,\wh{Q}_r \in R$ 
under the natural surjection 
$R \ra F$. 
Then,  
$(K,R,F)$ 
turns out to be a modular system. 
We denote by 
$\He_{n,r}^R$ 
(resp. $\Sc_{n,r}^R$) 
the   Ariki-Koike algebra 
(resp. cyclotomic $q$-Schur algebra) over 
$R$ 
with the parameters 
$\wh{q},\wh{Q}_1, \cdots, \wh{Q}_r$. 
Put 
$\He_{n,r}^{K}= K \otimes_R \He_{n,r}^R$ (resp. $\Sc_{n,r}^K= K \otimes_R \Sc_{n,r}^R$) 
and 
$\He_{n,r}^{F}= F \otimes_R \He_{n,r}^R$ (resp. $\Sc_{n,r}^F= F \otimes_R \Sc_{n,r}^R$), 
where we regard $F$ as the $R$-module through the natural surjection $R \ra F$. 
By using the modular system, 
for the Weyl module $W^\la$ ($\la \in \vL_{n,r}^+$) 
of $\Sc_{n,r}^F$, 
we can define the Jantzen filtration 
\[ 
W^\la = W^\la(0) \supseteq W^\la(1) \supseteq W^\la(2) \supseteq \cdots, 
\]
where 
$W^\la(1) = \rad W^\la$. 
Similarly, 
we have the Jantzen filtration 
\[ 
S^\la = S^\la(0) \supseteq S^\la(1) \supseteq S^\la(2) \supseteq \cdots
\]
for  
the Specht module $S^\la$ of $\He_{n,r}^F$.  

For 
$x=(i,j,k) \in [\la]$ ($\la \in \vL_{n,r}^+$), 
put 
\[ 
\res_{R} (x) = \wh{q}^{\,j-i}\,  \wh{Q}_k.
\]

For 
$\la,\mu \in \vL_{n,r}^+$, 
we define the integer $J_{\la\mu}$ 
called a \textbf{Jantzen coefficient} 
by 
\begin{align}
\label{def Jantzen coeff}
J_{\la\mu} = 
\begin{cases} 
\displaystyle 
\sum_{x \in [\la]} 
\sum_{y \in [\mu] \atop [\mu] \setminus r_y = [\la] \setminus r_x} 
(-1)^{\ell \ell (r_x) + \ell \ell(r_y)} 
\nu_{\wp} (\res_{R}(f_x) - \res_{R}(f_y)), 
& \text{if } \la \triangleright \mu, 
\\
0, 
& \text{otherwise,} 
\end{cases}
\end{align}
where 
$r_x$ is the rim hook of $x$, 
$\ell \ell (r_x)$ is the leg length of $r_x$ 
and 
$f_x$ is the foot node of $r_x$ 
(for definitions, see \cite{JM00}). 
Here, we remark that 
$f_x$ is a node in $\la^{(k)}$ 
if 
$x$ is a node in $\la^{(k)}$.   
Then we have the following theorem. 

\begin{thm}[{\cite{JM00}}] 
\label{thm Jantzen sum}
Let 
$(K,R,F)$ 
be a modular system. 
Suppose that 
$\Sc_{n,r}^K$ 
is semisimple.  
Then,  
for $\la \in \vL_{n,r}^+$,   
 we have 
\[ 
\sum_{i > 0} [ W^\la(i)] = \sum_{\mu \in \vL_{n,r}^+ \atop \la \triangleright \mu} 
J_{\la\mu} [W^\mu] 
\qquad (resp. 
\sum_{i > 0} [ S^\la(i)] = \sum_{\mu \in \vL_{n,r}^+ \atop \la \triangleright \mu} 
J_{\la\mu} [S^\mu]). 
\]
in the Grothendieck group of $\Sc_{n,r}^F$ (resp. $\He_{n,r}^F$).  
\end{thm}

\para
Let 
$e \in \{1,2,\cdots,\infty\}$ 
be the multiplicative order of $q$ in $F$, 
namely 
$q$ is a primitive $e$-th root of unity 
or 
$e=\infty$ 
if $q$ is not a root of unity. 
In the case where  
$r=1$, 
$\He_{n,1}(q,Q_1)$ 
is the Iwahori-Hecke algebra 
$\He_n(q)$ of type A, 
and 
$\Sc_{n,1}(q,Q_1)$ 
is the $q$-Schur algebra 
$\Sc_{n}(q)$ 
associated to $\He_n(q)$. 
(Note that 
$\He_{n,1}(q,Q_1)$ 
(resp. 
$\Sc_{n,1}(q,Q_1)$) 
is independent from the parameter 
$Q_1$.)
A condition for 
$\He_n(q)$ 
to be of finite representation type 
was obtained by Uno \cite{Uno92}. 
On the other hand, 
the representation type of $\Sc_n(q)$ 
has been determined by Erdmann-Nakano  
\cite{EN01} as follows.  

\begin{thm}[{\cite{Uno92}, \cite{EN01}}]
\label{Theorem EN} 
Suppose that $q\not=1$  
and $r=1$, 
then we have the following. 

\begin{enumerate}

\item 
$\He_{n}(q)$ (resp. $\Sc_{n}(q)$)  
is semisimple 
if and only if 
$n<e$. 

\item
$\He_{n}(q)$ (resp. $\Sc_{n}(q)$)  
is of finite type 
if and only if 
$n<2e$. 

\item
$\He_{n}(q)$ (resp. $\Sc_{n}(q)$)  
is of wild type 
if and only if 
$n\geq 2e$. 
\end{enumerate}
\end{thm}

\remarks\\
(\roi)\, 
In \cite{EN01A}, 
the representation type for an each block of $\He_n(q)$ 
is determined. 
\vspace{2mm}\\
(\roii)\,  
In this paper, 
we are only concerned with  cyclotomic $q$-Schur algebras satisfying the condition 
(CL)  
in \ref{comb}. 
In \cite{EN01}, 
the representation type of $q$-Schur algebras with 
$m_1 <n$  
is also determined. 
It occurs that 
$\Sc_{n}(q)$ 
is of tame type 
for some cases with 
$m_1 < n$. 
But, under the condition (CL), 
no $\Sc_n(q)$ has tame type. \vspace{2mm}\\
(\roiii)\, 
In the case where $q=1$, 
$q$-Schur algebras are nothing but classical Schur algebras. 
The representation type of (classical) Schur algebras 
has been determined by 
Doty-Erdmann-Martin-Nakano \cite{DEMN99}.

\para 
In order to describe the representation type 
of Ariki-Koike algebras and 
cyclotomic $q$-Schur algebras 
in the case where 
$r \geq 2$, 
we need the following theorem 
which has been proved by Dipper-Mathas 
\cite{DM02}.

\begin{thm}[{\cite{DM02}}]
\label{Th of DM}
Suppose that 
$I=I_1 \cup I_2 \cup \cdots \cup I_{\kappa} $ 
(disjoint union) 
is a partitioning of the index set 
$I=\{1,\cdots,r\}$ of parameters $\BQ=(Q_1,\cdots,Q_r)$ 
such that  
\[
\prod_{1\leq \a < \b \leq \kappa} 
\prod_{Q_i \in \BQ_\a \atop Q_j \in \BQ_\b}
\prod_{-n<a<n}
(q^a Q_i-Q_j) \not=0,\]
where 
we put 
$\BQ_\a=(Q_{\a_1},\cdots, Q_{\a_{j}})$ for $I_\a=\{\a_1,\cdots,\a_{j}\}$. 
Then 
$\He_{n,r}(q,\BQ)$ 
(resp. $\Sc_{n,r}(q,\BQ) $ ) 
is Morita equivalent to the algebra 
\[
\bigoplus_{n_1,\cdots,n_\kappa \geq 0 \atop n_1+\cdots +n_\kappa=n}
\He_{n_1,r_1}(q,\BQ_1)\otimes \He_{n_2,r_2}(q,\BQ_2) \otimes \cdots \otimes \He_{n_\kappa,r_\kappa}(q,\BQ_\kappa),
\]

\[\left( \text{resp. } 
\bigoplus_{n_1,\cdots,n_\kappa \geq 0 \atop n_1+\cdots +n_\kappa=n}
\Sc_{n_1,r_1}(q,\BQ_1)\otimes \Sc_{n_2,r_2}(q,\BQ_2) \otimes \cdots \otimes \Sc_{n_\kappa,r_\kappa}(q,\BQ_\kappa) 
\right),
\]
where $r_i=|\BQ_i|$ $(i=1,\cdots,\kappa)$. 
\end{thm}

\para 
By 
Theorem \ref{Th of DM}  
combined with 
Lemma \ref{tensor of cellular} 
(together with  multiplying $Q_1,\cdots, \\ Q_r$  by a scalar 
$c\in F$ 
simultaneously 
(see Remark \ref{remark} (\roiii))), 
we may assume that 
$Q_i=q^{f_i}$ ($i=1,\cdots,r$) or 
$Q_1=\cdots =Q_r=0$. 
Then, 
in order to determine the representation type of 
cyclotomic $q$-Schur algebras, 
it is enough to consider the following cases. 

\begin{description}

\item[case 1]
$q\not=1$ 
and 
$Q_i=q^{f_i}$ 
($f_i \in \ZZ$) 
for 
$i=1,\cdots,r$. 

\item[case 2]
$q=1$ and $Q_1=\cdots=Q_r=1$. 

\item[case 3]
$q=1$ and $Q_1=\cdots=Q_r=0$. 

\item[case 4]
$q\not=1$ and $Q_1=\cdots=Q_r=0$. 
\end{description}
In this paper, 
we are only concerned with  case 1 and  case 2.  
First, 
we consider  case 2.

\begin{thm}
Suppose that 
$q=1$, $r\geq 2$  
and 
$Q_1=\cdots=Q_r=1$. 
Then 
$\He_{n,r}(q,Q_1,\cdots,Q_r)$ 
(resp. 
$\Sc_{n,r}(q,Q_1,\cdots,Q_r)$) 
is of finite type 
if and only if 
$n=1$. 
\end{thm}

\begin{proof}
It is clear that 
$\He_{n,r}$ 
(resp. $\Sc_{n,r}$)  
is of finite type 
if
$n=1$. 
Suppose that $n\geq 2$. 
Then, 
by 
\cite[Proposition 41]{Ari05}, 
$\He_{n,2}(q,Q_1,Q_2)$ ($q=Q_1=Q_2=1$) 
is of  infinite type.
(Note that 
the parameters are given by 
$Q_1=Q_2=-1$ 
in \cite{Ari05}. 
By multiplying  $Q_1$ and $Q_2$ by $-1$, 
we obtain the above claim.) 
Thus 
$\He_{n,r}(q,Q_1,\cdots,Q_r)$ 
is of infinite type by Corollary \ref{cor induction}, 
and 
$\Sc_{n,r}(q,Q_1,\cdots,Q_r)$ 
is of infinite type by Lemma \ref{lem Hecke Schur}.
\end{proof}

\remark 
We can also prove that 
$\He_{n,r}$ 
(resp. $\Sc_{n,r}$)  
is of wild type 
if $r\geq 2$, $n\geq 3$ and $q=Q_1= \cdots =Q_r=1$ 
by \cite[Proposition 41]{Ari05}
in a similar way as in the above proof. 
But, we don't know 
whether 
$\He_{n,r}$ 
(resp. $\Sc_{n,r}$)  
is of tame type or of wild type 
if $n=2$, $r\geq 2$ and $q=Q_1=\cdots=Q_r=1$.

\para
From now on, we concentrate on  case 1 
with 
$r\geq 2$. 
Hence we assume the following condition. 
\begin{description}

\item[(CP)]
$q$ is a primitive $e$-th root of unity. 
($e=\infty$ if $q$ is not a root of unity. )

\quad  $Q_i=q^{f_i}$ $(0 \leq f_i \leq e-1)$ for $i=1,\cdots,r$. 
\end{description}
Note that, when $e=\infty$,  
we can take $f_i \in \ZZ_{\geq 0}$ 
without loss of generality by Remark \ref{remark} (\roiii), 
and  we regard as 
$c < \infty$ 
for any integer $c$. 
Let 
\[ 0 \leq f_1' \leq f_2' \leq \cdots \leq f_r'  \leq e-1 \]
be the increasing sequence of integers 
such that 
$f_i' = f_{\sigma(i)}$ ($i=1,\cdots,r$)  
for some permutation $\sigma$ of $r$ letters. 
Set $f'_{r+i}=e+f'_i$ 
and 
$g'_i=f'_{i+1}-f'_i$ 
($i=1,\cdots,r$). 
We define the integers 
$f^{+1}(Q_1,\cdots,Q_r)$, 
$f^{+2}(Q_1,\cdots,Q_r)$ 
and
$g(Q_1,\cdots,Q_r)$  
for parameters 
$Q_1,\cdots,Q_r$ 
by 
\begin{align*}
&f^{+1}(Q_1,\cdots,Q_r)= 
\min \big\{f'_{i+1}-f'_i\,|\,i=1,\cdots,r \big\}, \\
&f^{+2}(Q_1,\cdots,Q_r)=
\min \big\{f'_{i+2}-f'_i\,|\,i=1,\cdots,r \big\}, \\
&g(Q_1,\cdots,Q_r)= 
\min \big\{g'_i+g'_{j}\bigm| 1\leq i\not=j \leq r \big\}. 
\end{align*}
The rest of this section is devoted to the proof of the following theorem.

\begin{thm}
\label{main thm}
Under the condition (CP), we have the following. 

\begin{enumerate}
\item 
Assume that 
$r=2$. 
Then 
$\He_{n,2}(q,Q_1,Q_2)$ 
(resp. $\Sc_{n,2}(q,Q_1,Q_2)$)  
is of finite type 
if and only if 
\[
n<\min \{e,\, 2f^{+1}(Q_1,Q_2)+4\}.
\]

\item 
Assume that 
$r\geq 3$. 
Then 
$\He_{n,r}(q,Q_1,\cdots,Q_r)$ 
(resp. $\Sc_{n,r}(q,Q_1,\cdots,Q_r)$)  
is of finite type 
if and only if 
\[
n<\min\{2f^{+1}(Q_1,\cdots,Q_r)+4,\, f^{+2}(Q_1,\cdots,Q_r)+1, \, g(Q_1,\cdots,Q_r)+2\}.
\]

\end{enumerate}
\end{thm}

%
%
%
%
%


\para
In order to prove Theorem \ref{main thm}, 
we prepare some known results on blocks of 
$\He_{n,r}$ and $\Sc_{n,r}$. 
By a general theory of cellular algebras 
\cite{GL96}, 
for each 
$\la \in \vL^+_{n,r}$, 
all of the composition factors of 
the Specht module 
$S^\la$ of $\He_{n,r}$ 
belong to the same block of 
$\He_{n,r}$. 
This result allows us to classify the blocks of 
$\He_{n,r}$ 
by using the Specht modules. 
Similar facts also holds for 
$\Sc_{n,r}$. 
By Lyle-Mathas \cite{LM07}, 
this classification 
has been  described 
combinatorially as follows.  
Here, 
we only give the result under the condition (CP) 
though 
it is described in a general setting 
in \cite{LM07}.

For 
$\la \in \vL_{n,r}^+$, 
we define the 
\textit{residue} 
of the node 
$x=(i,j,k)\in [\la]$ 
by 
\[
\res (x) =j-i+f_k \quad(\hspace{-1em}\mod e). 
\]
For 
$\la,\mu \in \vL_{n,r}^+$, 
we say that 
$\la$ 
and  
$\mu$ 
are \textit{residue equivalent}, 
and write 
$\la \sim_C \mu$ 
if 
$\sharp \{x \in [\la] \,|\,res(x)=a\} 
= \sharp \{ y \in [\mu] \,|\, \res(y)=a\}$ 
for all $a \in \ZZ/e \ZZ$. 

\begin{thm}[{\cite[Theorem 2.11]{LM07}}]
\label{Thm block}
For 
$\la,\mu \in \vL_{n,r}^+$, 
the following conditions are equivalent. 

\begin{enumerate}

\item
$S^\la$ and $S^\mu$ 
belong to the same block of 
$\He_{n,r}(q,Q_1,\cdots,Q_r)$. 

\item
$W^\la$ and $W^\mu$ 
belong to the same block of 
$\Sc_{n,r}(q,Q_1,\cdots,Q_r)$. 

\item
$\la \sim_C \mu$. 
\end{enumerate}
\end{thm}

For a block 
$\CB$ 
of $\Sc_{n,r}$, 
we write 
$\la \in \CB$ 
if 
$W^\la$ 
belongs to the block 
$\CB$, 
and we say that 
$\CB$ 
has 
a residue 
$(\res(x_1),\cdots,\res(x_m))$, 
where 
$\{x_i \,|\, i=1,\cdots, m\}=[\la]$ 
for some 
$\la \in \CB$. 
By Theorem \ref{Thm block}, 
the residue 
$(\res(x_1),\cdots,\res(x_m))$ 
of $\CB$ 
is well-defined 
up to a permutation of components. 
Set 
$R(\CB)=\{\res(x)\,|\,x \in [\la]  
\text{ for some } 
\la \in \CB \}$. 
It is similar for $\He_{n,r}$, 
and we use the same notations. 

Recall that 
$F$ is the Schur functor defined in \ref{Schur functor}.
As a corollary of Theorem \ref{Thm block}, 
we have the following. 

\begin{cor}
\label{cor Schur functor}
For $\la \in \vL^+_{n,r}$, 
let $P^\la$ be the projective cover of $L^\la$. 
Then we have the following. 

\begin{enumerate}
\item 
$F(P^\la) \not=0$ for any $\la \in \vL_{n,r}^+$. 

\item 
If 
$F(L^\la)\not=0$ and $F(L^\mu) \not=0$, 
then the following two conditions are equivalent. 
\begin{enumerate}
\item 
$P^\la$ and $P^\mu$ belong to the same block of $\Sc_{n,r}$. 
\item
$F(P^\la)$ and $F(P^\mu)$ belong to the same block of $\He_{n,r}$. 
\end{enumerate}
\end{enumerate}
\end{cor}

\begin{proof}
By the general theory of quasi-hereditary algebras, 
there exists a submodule $K^\la$ of $P^\la$ 
such that 
$P^\la / K^\la \cong W^\la$ 
for each $\la \in \vL^+_{n,r}$. 
Since $F$ is an exact functor, 
we have 
that 
$F(P^\la)/ F(K^\la) \cong F(W^\la) \cong S^\la$ 
by Lemma \ref{lem Schur functor} (\roi). 
This implies (\roi), 
and 
(\roii)  
follows from 
Theorem \ref{Thm block}. 
\end{proof}
\para 
First, 
we prove the \lq\lq \,if " part of  
Theorem \ref{main thm}. 
The following Lemma 
plays an important role 
in the proof of the \lq\lq \, if " part of 
Theorem \ref{main thm}. 

\begin{lem}
\label{Lemma fin}
Under the condition (CP), 
suppose that 
$r\geq 3$ 
and 
\begin{align}
\label{n < f-g}
n< \min \big\{f^{+2}(Q_1,\cdots,Q_r)+1, \, g(Q_1,\cdots,Q_r)+2\big\}.
\end{align}
Then, for each block 
$\CB$ 
of 
$\Sc_{n,r}(q,Q_1,\cdots,Q_r)$,  
there exist  
$i,j \in \{1,\cdots,r\}$ 
such that 
$\la^{(k)}=\mu^{(k)}$ 
for any 
$\la,\mu \in \CB$, 
and for  
$k\not= i,j$.
\end{lem}

\begin{proof}
Let $\sigma$ 
be a permutation of $r$ letters. 
We define the bijective map $\sigma : \vL_{n,r}^+ \ra \vL_{n,r}^+$ 
by 
$\s((\la^{(1)},\la^{(2)},\cdots,\la^{(r)})) = (\la^{(\s^{-1}(1))},\la^{(\s^{-1}(2))}, \cdots, \la^{(\s^{-1}(r))}) $. 
For $\la \in \vL_{n,r}^+$, 
we also define the bijection 
$\s : [\la] \ra [\s(\la)] $ by 
$\s(x)=(i,j,\s(k))$  
for $x=(i,j,k) \in [\la]$. 
Then, 
for $\la \in \vL_{n,r}^+$ and $x=(i,j,k) \in [\la]$, 
one sees easily that 
$\res (x)$ with respect to the parameters 
$q,Q_1,\cdots,Q_r$ 
coincides with 
$\res (\s(x))$ 
with respect to the parameters 
$q,Q_{\s(1)},\cdots, Q_{\s(r)}$. 
Combining with Theorem \ref{Thm block}, 
we have that 
$\la,\mu$ belong to the same block of 
$\Sc_{n,r}(q,Q_1,\cdots,Q_r)$ 
if and only if 
$\la,\mu$ belong to the same block of 
$\Sc_{n,r}(q,Q_{\s(1)}, \cdots, Q_{\s(r)})$. 
Hence, 
it is enough to show the following cases:  
\begin{align} 
\label{condition f}
0 \leq f_1 \leq f_2 \leq \cdots \leq f_r \leq e-1.
\end{align}
Hence, 
in this proof, 
we suppose that the condition 
\eqref{condition f} holds. 
Note that 
$f_i'=f_i$ ($i=1,\cdots,r$) 
under the condition \eqref{condition f}.

Take $\la, \mu \in \CB$ such that $\la\not=\mu$. 
For a node 
$x=(a,b,i) \in [\la] \setminus [\mu]$, 
there exists a node 
$y=(c,d,j) \in [\mu] \setminus [\la]$ 
such that 
$\res(x)=\res(y)$
by 
Theorem \ref{Thm block}. 
Since 
the condition 
$n < f^{+2}(Q_1,\cdots,Q_r) +1$ 
implies that 
$n < e$ by a direct calculation, 
we have 
$i \not=j$. 
We may assume that 
$i<j$ 
by 
interchanging  
$\la $ and $\mu$ 
if necessary. 
Since 
$\res(x)=\res(y)$, 
one of the following two cases occurs:   
\begin{align}
\label{RB contain}
\begin{cases}  
R(\CB) \supseteq \{f_{i},f_{i}+1 ,\cdots,f_{j}\},\\ 
R(\CB) \supseteq \{f_{j}, f_{j}+1, \cdots, e-1, 0, 1, \cdots, f_{i}\}.\\
\end{cases}
\end{align}
In fact, 
suppose that 
$a \geq b$, 
then 
$f_i, f_{i}+1, \cdots, \res(x)$ 
occur among  the residues of 
$\la$. 
If 
$c \geq d$, 
then 
$f_{j}, f_{j}+1, \cdots, \res(x)$ 
occur in 
$\mu$, 
and 
if 
$c<d$, 
then 
$\res(x),\res(x)+\nobreak1,\cdots, f_j$ 
occur in 
$\mu$.   
In either case, 
we see that 
$f_i,\cdots, f_j$ 
occur in 
$R(\CB)$. 
Next suppose 
$a<b$, 
then 
$\res(x), \res(x)+1,\cdots, f_i$, 
occur in 
$\la$. 
Hence 
if 
$c \geq d$, 
again 
$f_i, f_i+1, \cdots, f_j$ 
occur in 
$R(\CB)$, 
and 
if 
$c<d$, 
we see that 
$\{f_j, \cdots, e-1, 0, 1, \cdots, f_i\} \in R(\CB)$. 

If 
$|i-j|>1$ and $\{i,j\}\not=\{1,r\}$,  
any case of (\ref{RB contain}) 
implies that 
\begin{align}
\label{n > f+2}
n \geq f^{+2}(Q_1,\cdots, Q_r)+1.
\end{align}
In fact, 
the first case implies that 
$f_j-f_i+1 \leq n$. 
In the second case, 
we have 
$e- f_j+f_i \leq n-1$, 
which implies that 
$n \geq f_{j+2}-f_{j}+1$ 
since
$\{i,j\} \not=\{1,r\}$. 
But 
(\ref{n > f+2}) contradicts the condition (\ref{n < f-g}). 
Thus we have 
\[
|i-j|=1 \text{ or } \{i,j\}=\{1,r\}.
\] 

Next 
we show that 
such a pair 
$\{i,j\}$ 
is determined uniquely 
for a given 
$\la,\mu \in \CB$. 
Let 
$\la, \mu \in \CB$, 
and  
take  
$x_k=(a_k,b_k,i_k) \in [\la] \setminus  [\mu] $, 
$y_k=(c_k,d_k,j_k) \in [\mu] \setminus [\la]$ 
such that 
$\res(x_k)=\res(y_k)$ 
($k=1,2$). 
By the above result, 
we have 
$|i_k-j_k|=1$ 
or 
$\{i_k,j_k\}=\{1,r\}$. 
If 
$\{i_1,j_1\}\cap \{i_2,j_2\}=\emptyset$, 
we have 
$n \geq g(Q_1,\cdots,Q_r)+2$ 
by considering the residues contained  in 
$R(\CB)$ 
of (\ref{RB contain}). 
This contradicts the condition (\ref{n < f-g}). 
For example, 
in the case where 
$r=6, i_1=1,j_1=2, i_2=3,j_2=4$, 
we have 
\[R(\CB)\supset \{f_1,f_1+1,\cdots,f_2 \} \cup \{f_3,f_3+1,\cdots,f_4 \}.\]
Thus we may assume that 
\begin{align}
\label{i,j-empty}
\{i_1,j_1\}\cap \{i_2,j_2\}\not=\emptyset. 
\end{align}
If 
the two sets 
contain exactly one common element,  
we have 
$n \geq f^{+2}(Q_1, \cdots,Q_r)+1$  
by considering the residues contained  in 
$R(\CB)$ of (\ref{RB contain}). 
This contradicts the condition (\ref{n < f-g}). 
For example, 
in the case where 
$r=6, i_1=1,j_1=2, i_2=2,j_2=3$, 
we have 
\[R(\CB)\supset \{f_1,f_1+1,\cdots,f_2 \} \cup \{f_2,f_2+1,\cdots,f_3 \}.\]
As a conclusion, 
we have 
$\{i_1,j_1\}=\{i_2,j_2\}$.

Finally, 
we show that 
such a pair 
$\{i,j\}$ 
is independent of the choice of 
$\la,\mu \in \CB$. 
For 
$\la,\mu,\nu \in \CB$, 
take 
$x=(a,b,i) \in [\la] \setminus [\mu]$, 
$y=(c,d,j) \in [\mu] \setminus [\la]$ 
such that 
$\res(x)=\res(y)$,  
and take 
$x'=(a',b',i') \in [\mu] \setminus [\nu]$, 
$y'=(c',d',j') \in [\nu] \setminus [\mu]$ 
such that 
$\res(x')=\res(y')$. 
If 
$\{i,j\}\cap \{i',j'\}=\emptyset$, 
we have 
$n \geq g(Q_1,\cdots,Q_r)+2$, 
and 
if  
$\{i,j\}$ 
and 
$\{i',j'\}$ 
contain exactly one common element, 
we have 
$n \geq f^{+2}(Q_1, \cdots,Q_r)+1$ 
in a similar way as above. 
Thus we have 
$\{i,j\}=\{i',j'\}$. 
The lemma is proved. 
\end{proof}


Lemma \ref{Lemma fin} implies the following proposition 
which shows the \lq\lq\, if " part of Theorem \ref{main thm}.

\begin{prop}
\label{finiteness}
Under the condition (CP), 
if 
$n<\min \{2f^{+1}(Q_1,\cdots,Q_r)+4, \, f^{+2}(Q_1,\cdots,Q_r)+1, \, g(Q_1,\cdots, Q_r)+2\}$ 
$(r\geq 3)$ 
or if 
$n<\min \{e,\, 2f^{+1}(Q_1,Q_2)+4\}$ 
$(r=2)$, 
then 
the decomposition matrix of 
a block of 
$\Sc_{n,r}(q,Q_1,\cdots,Q_r)$ 
is given as 
\[D=\left(
\begin{matrix}
1&&&\\
1&1&&\\
&\hspace{-1mm}\ddots\hspace{-1mm}&\hspace{-1mm}\ddots\hspace{-1mm}&\\
&&1&1
\end{matrix}
\right)
\quad
\text{(all omitted entries are zero)}. 
\]
Moreover, 
any projective indecomposable $\Sc_{n,r}$-module 
has the simple socle. 

In particular, 
$\Sc_{n,r}(q,Q_1,\cdots,Q_r)$ 
is of 
finite type. 
Thus, 
$\He_{n,r} (q,Q_1,\cdots,Q_r)$ 
is also of 
finite type. 
\end{prop}

\begin{proof}
Take  a block 
$\CB$ 
of 
$\Sc_{n,r}(q,Q_1,\cdots,Q_r)$. 
We compute the decomposition matrix 
of 
$\CB$ 
by 
the Jantzen sum formula 
in Theorem \ref{thm Jantzen sum}. 

In the case where  
$r\geq3$, 
assume that 
$n < \min \{ 2 f^{+1}(Q_1,\cdots,Q_r)+4, f^{+2}(Q_1,\\   \cdots,Q_r)+1, \, g(Q_1,\cdots,Q_r)+2\}$. 
Then by Lemma \ref{Lemma fin}, 
there exist  
$i,j \in \{1,\cdots,r\}$ 
such that, 
$\la^{(k)}=\mu^{(k)}$ 
for any 
$\la,\mu \in \CB$ 
and 
$k\not= i,j$.
Put 
$n'=|\la^{(i)}|+|\la^{(j)}|$ 
for $\la \in \CB$, 
which is independent of 
$\la \in \CB$. 
By Theorem \ref{Thm block} combined with Lemma \ref{Lemma fin}, 
one can find a block 
$\CB'$ 
of 
$\Sc_{n',2}(q,Q_i,Q_{j})$ 
which contains 
$\{(\la^{(i)},\la^{(j)})\,|\,\la \in \CB\}$. 
In order to compute the decomposition matrix of $\CB$ by the Jantzen sum formula, 
we take the following modular system.  
Let 
$F[t]$ 
be a polynomial ring over $
F$ 
with indeterminate 
$t$, 
and 
$R=F[t]_{\lan t \ran}$ 
be the localization of 
$F[t]$ 
by the prime ideal 
$\lan t \ran$ 
generated by 
the polynomial $t$.  
Let 
$K$ 
be 
the quotient field of 
$R$. 
Put 
$\wh{q}=q$, 
$\wh{Q}_i= q^{f_i}$, 
$\wh{Q}_j = t+ q^{f_j}$ 
and 
$\wh{Q}_k = t^{2k} + q^{f_k}$ for $k\not=i,j$  
as elements in $R$. 
Then,  
$(K,R,F)$ 
becomes a modular system. 
Under this modular system, 
 we see that $\Sc_{n,r}^K(\wh{q},\wh{Q}_1,\cdots, \wh{Q}_r)$ 
is semisimple. 
Thus, we can apply the Jantzen sum formula (Theorem \ref{thm Jantzen sum}). 
By the definition \eqref{def Jantzen coeff} 
combined with Lemma \ref{Lemma fin}, 
the Janzten coefficient 
$J_{\la\mu}$ 
for 
$\la,\mu \in \CB$ 
is determined by only the informations of 
$(\la^{(i)},\la^{(j)})$ 
and 
$(\mu^{(i)},\mu^{(j)})$ 
since 
$\la^{(k)}=\mu^{(k)}$ 
for 
$k\not= i,j$. 
Moreover, 
this Jantzen coefficient 
$J_{\la\mu}$ 
coincides with 
the Jantzen coefficient 
$J_{(\la^{(i)},\la^{(j)}),(\mu^{(i)},\mu^{(j)})}$
in 
$\Sc_{n',2}(q,Q_i,Q_{j})$, 
where we take the same modular system 
$(K,R,F)$ 
with the parameters 
$\wh{q},\wh{Q}_i, \wh{Q}_j$. 
This means that 
the Jantzen sum formula for 
$\CB$ 
coincides with 
the Jantzen sum formula for 
$\CB'$. 
Moreover, 
$\Sc_{n',2}(q,Q_i,Q_{j})$ 
satisfies the assumption of the proposition 
for the case where $r=2$. 
Thus, 
we have only to compute the decomposition matrix in the case where  
$r=2$. 

Note that, 
the Jantzen sum formula 
(more precisely the Jantzen coefficient)  
for 
$\Sc_{n,r}(q,Q_1,\cdots,Q_r)$ 
coincides with 
the Jantzen sum formula for 
$\He_{n,r}(q,Q_1,\cdots,Q_r)$ 
(see Theorem \ref{thm Jantzen sum}). 
For $\He_{n,2}(q,Q_1,Q_2)$ 
satisfying the condition 
$n<\min \{e,\, 2f^{+1}(Q_1, \\ Q_2)+4\}$, 
the Jantzen coefficient has been computed by 
\cite[Theorem 6.2]{AM04}. 
Thus, 
we can compute the decomposition matrix of 
$\CB$ 
in a similar way as in the proof of [loc. cit.], 
and we obtain the matrix as given in the proposition. 

Next, 
we show that 
any projective indecomposable $\Sc_{n,r}$-module 
has the simple socle. 
Let 
$\{\la_1, \cdots, \la_m\} =\{\la \in \CB\}$ 
be such that $i<j$ if $\la_i \trl \la_j$. 
Then 
$P^{\la_i}$ has the radical series as in \eqref{radical layer}. 
It is clear that 
$P^{\la_m}$ 
has the simple socle 
from the radical series. 
%
%
We claim that 
\begin{align}
\label{claim F not vanish} 
F(L^{\la_{i}}) \not=0 \text{ for } i=1,\cdots,m-1.
\end{align} 
By Lemma \ref{lem Schur functor} (\roi), 
we see that 
$F(L^{\la_{1}})\not=0$ 
since 
$L^{\la_1}=W^{\la_1}$. 
Assume that $F(L^{\la_i})=0$ for some $i=2,\cdots,m-1$. 
Then  
we have that 
$F(L^{\la_{i-1}}) \not=0$ or $F(L^{\la_{i+1}})\not=0$ 
since $F(P^{\la_i})\not=0$ by Corollary \ref{cor Schur functor} (\roi). 
If 
$F(L^{\la_{i-1}}) \not=0$ and $F(L^{\la_{i+1}}) \not=0$, 
then  
$F(P^{\la_{i-1}})$ and $F(P^{\la_{i+1}})$ 
are projective indecomposable $\He_{n,r}$-modules  
by Lemma \ref{lem Schur functor} (\roiii). 
Moreover, 
$F(P^{\la_{i-1}})$ and $F(P^{\la_{i+1}})$ 
belong to the same block of $\He_{n,r}$ 
by Corollary \ref{cor Schur functor} (\roii). 
However, 
$F(P^{\la_{i-1}})$ and $F(P^{\la_{i+1}})$ 
are not linked 
since $F(L^{\la_i})=0$. 
This is a contradiction.  
If  
$F(L^{\la_{i-1}}) =0$, 
we have that 
$F(L^{\la_{i+1}})\not=0$, 
and also  
$F(L^{\la_{i-2}})\not=0$ 
since $F(P^{\la^{i-1}})\not=0$. 
Then $F(P^{\la_{i-2}})$ and $F(P^{\la_{i+1}})$ 
belong to the same block of $\He_{n,r}$, 
but these are not linked. 
This is a contradiction.
In the case where $F(L^{\la_{i+1}})=0$, 
we have a similar contradiction. 
Hence, 
we have the claim \eqref{claim F not vanish}. 

By \eqref{claim F not vanish} combined with Lemma \ref{lem Schur functor} (\roiii), 
$F(P^{\la_{i}})$ $(1 \leq i \leq m-1)$ 
is a projective indecomposable $\He_{n,r}$-module,  
and it is also injective 
since $\He_{n,r}$ is self-injective by \cite{MM}. 
Thus, $F(P^{\la_{i}})$ has the simple socle. 
Combining with \eqref{claim F not vanish},  
we see that 
$P^{\la_i}$ ($1 \leq i \leq m-1$) has the simple socle. 


Finally, 
we note that the block 
$\CB$ 
of 
$\Sc_{n,r}(q,Q_1,\cdots,Q_r)$ 
is also a quasi-hereditary cellular algebra.    
Then  
we conclude that 
$\Sc_{n,r}(q,Q_1,\cdots, Q_r)$ 
is of finite type 
by Proposition \ref{decom matrix}, 
and 
$\He_{n,r}(q,Q_1,\cdots,Q_r)$ 
is also of finite type 
by Lemma \ref{lem Hecke Schur}.
\end{proof}

As a corollary of the proof of the proposition, 
we have the following. 

\begin{cor}
\label{Morita}
Suppose that $r \geq 3$,  
and that 
$\Sc_{n,r}(q,Q_1,\cdots,Q_r)$ 
satisfies the assumption of 
Proposition $\ref{finiteness}$. 
Then, 
for each block 
$\CB$ 
of 
$\Sc_{n,r}(q,Q_1,\cdots,Q_r)$, 
there exist $i,j \in \{1,\cdots,r\}$ 
such that 
$\CB$ 
is Morita equivalent to 
a block $\CB'$ of 
$\Sc_{n',2}(q,Q_i,Q_{j})$, 
where 
$n'=|\la^{(i)}|+|\la^{(j)}|$ 
for some $\la \in \CB$, 
and 
$\CB'$ contains the partitions 
$\{(\la^{(i)},\la^{(j)})\,|\,\la \in \CB\}$.  

\end{cor}

\begin{proof}
By the proof of Proposition \ref{finiteness} 
and 
Proposition \ref{decom matrix}, 
both of 
$\CB$ 
and 
$\CB'$ 
are Morita equivalent to 
the algebra $\CA_m$ in \ref{def Am}, 
where $m=\sharp\{\la \in \CB\}$. 
\end{proof}



Next, 
we discuss the 
\lq\lq only if " part of Theorem \ref{main thm}. 

\begin{prop}
\label{infinite S3}
Under the condition (CP),  
suppose that 
$n< \min\{e,\,2f^{+1}(Q_1, \\ Q_2, Q_3)+4\}$ 
and 
$n\geq f^{+2}(Q_1,Q_2,  Q_3)+1$. 
Then 
$\He_{n,3}(q,Q_1,Q_2, Q_3)$ 
is of infinite type. 
\end{prop}

\begin{proof}
By Remarks \ref{remark} (\roiii), 
we may assume that 
$0 = f_1 \leq f_2 \leq f_3 \leq e-1$ 
without loss of generality. 
First, 
we assume that 
$n \geq 4$. 
Then, 
by the condition 
$n<2f^{+1}(Q_1,Q_2,Q_3)+4$, 
we have 
$0=f_1 < f_2 < f_3$. 
Moreover, 
we can take 
$f^{+2}(Q_1,Q_2,Q_3)=f_3-f_1=f_3$ 
by permuting and multiplying $Q_1,Q_2,Q_3$ by a common  scalar  
if necessary. 
Let 
$\CB$ 
be the block of 
$\He_{n,3}(q,Q_1,Q_2,Q_3)$ 
with the residue 
$(f_3,f_3-1,\cdots,1,0,e-1,\cdots,e-f')$, 
where 
$f'=n-(f_3+1)$. 
The condition 
$n \geq f^{+2}(Q_1,Q_2,Q_3)+1=f_3+1$ 
implies that 
$f' \geq 0$, 
and 
the condition 
$n<e$ 
implies that 
$e-f'>f_3+1$. 

Put 
$\la_0=(-,-,(1^{n})),\,
\la_i=(-,(i,1^{f_2+f'}),(1^{f_3-f_2-i+1}))$ 
for $i=1,\cdots, k-1$ 
and 
$\la_{k}=((1^{n-f_3}),-,(1^{f_3})) $, 
where
$k=f_3-f_2+2$. 

\textbf{\underline{Example}} \,
(In the case where $n=5, \, e=6, \, f_1=0,\, f_2=1,\, f_3=3$.) \\

\hspace{3em}
$\la_0=(-,\,-,\,
{\tiny 
\begin{tabular}{|c|c|c|c|}
\hline 
3 \\ \hline 
2 \\ \hline 
1 \\ \hline 
0 \\ \hline
5 \\ \hline
\end{tabular}
}
\,\,)
$, \,
$\la_1=(-,\, 
{\tiny 
\begin{tabular}{|c|c|c|c|}
\hline 
1 \\ \hline 
0 \\ \hline 
5 \\ \hline 
\end{tabular}
}
\,,\,
{\tiny
\begin{tabular}{|c|c|c|c|}
\hline 
3 \\ \hline 
2 \\ \hline 
\end{tabular}
}
\,\,)
$,\,
$\la_2=(-,\, 
{\tiny 
\begin{tabular}{|c|c|c|c|}
\hline 
1 &2 \\ \hline 
0 \\ \cline{1-1} 
5 \\ \cline{1-1}
\end{tabular}
}
\,,\,
{\tiny
\begin{tabular}{|c|c|c|c|}
\hline 
3 \\ \hline
\end{tabular}
}
\,\,)
$,\\

\hspace{3em}
$\la_3=(-,\, 
{\tiny 
\begin{tabular}{|c|c|c|c|}
\hline 
1 &2 &3\\ \hline 
0 \\ \cline{1-1}
5 \\ \cline{1-1}
\end{tabular}
}\,,\,
-)
$,\,
$\la_4=(\,
{\tiny 
\begin{tabular}{|c|c|c|c|}
\hline 
0 \\ \hline
5 \\ \hline
\end{tabular}
}
\,, 
-,\,
{\tiny
\begin{tabular}{|c|c|c|c|}
\hline 
3 \\ \hline 
2 \\ \hline 
1 \\ \hline
\end{tabular}
}
\,\,)
$. \vspace{3mm}\\

It is easy to see that 
$\la_0, \cdots, \la_k$ all lie in the  block $\CB$ 
, and 
let $\la_{k+1},\cdots, \la_m$ 
denote the remaining multipartitions in $\CB$, in any order.  
Then 
$\CB$ 
is a cellular algebra 
with respect to the poset 
$(\vL_{\CB}^+, \trreq)$, 
where 
$\vL_{\CB}^+=\{\la_i \,|\, 0 \leq i \leq m\}$.   
Let 
$\CB^\vee$ 
be the $F$-subspace of 
$\CB$ 
spanned by the cellular basis elements 
indexed by 
$\la_i$ ($i\not= 0,1,k$). 
Note that 
$e-f'>f_3+1$ 
and 
the definitions of 
multipartitions 
$\la_0,\la_1,\cdots,\la_k$,  
one sees that 
$\la_i \ntriangleright \la_j$ 
for $i=0,1,k$ and $ j\not= 0,1,k$. 
This implies that 
$\CB^\vee$ 
is a two-sided ideal of 
$\CB$. 
Thus  
$\ol{\CB}=\CB/\CB^\vee$ 
becomes a cellular algebra 
with respect to the poset 
$(\{\la_0,\la_1,\la_k\}, \trreq)$. 
One can easily check that  
$\la_i$ ($i=0,1,k$) 
is a Kleshchev multipartition, 
thus 
$D^{\la_i} \not=0$ for $i=0,1,k$. 
Now, we can compute the decomposition matrix of 
$\ol{\CB}$ 
by the Jantzen sum formula, 
and its matrix is  given as follows.\\

\begin{tabular}{c|ccc}

&$D^{\la_0}$ &$D^{\la_1}$ & $D^{\la_k}$ \\ \hline 
$S^{\la_0}$&1&0&0\\
$S^{\la_1}$&a&1&0\\
$S^{\la_k}$&b&0&1\\
\end{tabular}
$\qquad (a,b >0)$
\\[3mm]
This implies that 
$\ol{\CB}$ is of infinite type 
by Lemma \ref{lem infinite cellular},   
thus 
$\CB$ 
and so 
$\He_{n,3}(q,Q_1,Q_2,Q_3)$ is of infinite type.

In the case where $n\leq 3$, 
we have the following three cases, 
$0=f_1=f_2=f_3$, $0=f_1=f_2 < f_3$ or $f_1 <  f_2 < f_3$. 
For each case, 
we can prove in a similar way 
as above  
by taking the appropriate block.  
\end{proof}

\begin{prop}
\label{infinite condition g}
Under the condition (CP),  
suppose that 
$n< \min\{2f^{+1}(Q_1,\cdots , \\ Q_4)+4, \, f^{+2}(Q_1,\cdots , Q_4)+1 \}$ 
and 
$n\geq g(Q_1, \cdots ,Q_4)+2$. 
Then 
$\He_{n,4}(q,Q_1,\cdots ,Q_4)$ 
is of infinite type. 
\end{prop}

\begin{proof}
By Remarks \ref{remark} (\roiii), 
we may assume that 
$0 = f_1 \leq f_2 \leq f_3 \leq e-1$ 
without loss of generality. 
First, we suppose that 
$n \geq 4$. 
Then, 
by the condition 
$n<2f^{+1}(Q_1,\cdots,Q_4)+4$, 
we have 
$0=f_1 < f_2 < f_3 < f_4$. 
If 
$g(Q_1,\cdots,Q_4)=g_i+g_{i+1}$ 
for some 
$i \in \{1,2,3,4\}$, 
we have 
$g(Q_1,\cdots,Q_4)=f^{+2}(Q_1,\cdots,Q_4)$. 
In this case, 
the condition 
$n \geq g(Q_1,\cdots, Q_4)+2$ 
contradicts 
the condition 
$n< f^{+2}(Q_1,\cdots,Q_4)+1$. 
Thus 
one can assume that 
$g(Q_1, \cdots, Q_4)=g_1+g_3=(f_2-f_1)+(f_4-f_3)$ 
(by permuting the parameters, 
and by a scalar multiplication 
if necessary).

Let 
$\CB$ 
be the block of 
$\He_{n,4}(q,Q_1,\cdots,Q_4)$ 
with the residue
$(0,\,1,\, 2,\cdots,f_2,\, f', \\   f'+1, \cdots,f_4-1,\, f_4)$, 
where 
$f'=f_4-(n-(f_2+1))+1$. 
Note that 
$(f_2+1)+(f_4-f'+1)=n$ 
and that 
$f'\leq f_3 <f_4$ 
by the condition 
$n \geq g_1+g_3+2$. 
Moreover, 
by the condition 
$n<f^{+2}(Q_1,\cdots,Q_r)+1$, 
we have 
$n<f_4-f_2+1 < f_4+1$.  
This implies that 
$f_2 < f'-1$. 
Put 
$\la_0=\big(-,(1^{f_2+1}),-,(1^{f_4-f'+1})\big)$, 
$\la_i=\big(-,(1^{f_2+1}),(i,1^{f_3-f'}), (1^{f_4-f_3-i+1})\big)$ 
for 
$i=1,\cdots, k$ 
and 
$\la_{k+1}=\big((1),(1^{f_2}),-,(1^{f_4-f'+1})\big)$, 
where 
$k=f_4-f_3+1$. 

\textbf{\underline{Example}} \,
(In the case where $n=7, \, e=16,\, f_1=0,\, f_2=2,\, f_3=8,\, f_4=10$.) \\

\hspace{1em}
$\la_0=(-,\,
{\tiny 
\begin{tabular}{|c|c|c|c|}
\hline 
2 \\ \hline 
1 \\ \hline 
0 \\ \hline 
\end{tabular}
}
,\,-,\,
{\tiny 
\begin{tabular}{|c|c|c|c|}
\hline 
10 \\ \hline 
9 \\ \hline 
8 \\ \hline 
7 \\ \hline
\end{tabular}
}
\,\,)
$, \,
$\la_1=(-,\, 
{\tiny 
\begin{tabular}{|c|c|c|c|}
\hline 
2 \\ \hline
1 \\ \hline 
0 \\ \hline 
\end{tabular}
}
\,,\,
{\tiny
\begin{tabular}{|c|c|c|c|}
\hline 
8 \\ \hline 
7 \\ \hline 
\end{tabular}
},
\, 
{\tiny
\begin{tabular}{|c|c|c|c|}
\hline 
10 \\ \hline 
9 \\ \hline 
\end{tabular}
}\,)
$,\,
$\la_2=(-,\, 
{\tiny 
\begin{tabular}{|c|c|c|c|}
\hline 
2 \\ \hline
1 \\ \hline 
0 \\ \hline 
\end{tabular}
}
\,,\,
{\tiny 
\begin{tabular}{|c|c|c|c|}
\hline 
8 &9 \\ \hline 
7 \\ \cline{1-1}
\end{tabular}
}
\,,\,
{\tiny
\begin{tabular}{|c|c|c|c|}
\hline 
10 \\ \hline
\end{tabular}
}
\,\,)
$,\\

\hspace{1em}
$\la_3=(-,\, 
{\tiny 
\begin{tabular}{|c|c|c|c|}
\hline 
2 \\ \hline
1 \\ \hline 
0 \\ \hline 
\end{tabular}
}
\,,\,
{\tiny 
\begin{tabular}{|c|c|c|c|}
\hline 
8 &9 &10\\ \hline 
7 \\ \cline{1-1}
\end{tabular}
}\,,\,
-)
$,\,
$\la_4=(\,
{\tiny 
\begin{tabular}{|c|c|c|c|}
\hline 
0 \\ \hline
\end{tabular}
}
\,, 
{\tiny 
\begin{tabular}{|c|c|c|c|}
\hline 
2 \\ \hline
1 \\ \hline 
\end{tabular}
}
\,,\,
-,\,
{\tiny
\begin{tabular}{|c|c|c|c|}
\hline 
10 \\ \hline 
9 \\ \hline 
8 \\ \hline
7 \\ \hline
\end{tabular}
}
\,\,)
$. \vspace{3mm}\\

In a similar way as in the proof of Proposition \ref{infinite S3}, 
one can consider the quotient algebra 
$\ol{\CB}=\CB/\CB^\vee$ 
of $\CB$ 
which is a cellular algebra  
with respect to the poset 
$\big(\{\la_0, \la_1, \la_{k+1}\},\,\trianglerighteq \big)$. 
One can easily check that  
$\la_i$ ($i=0,1,k+1$) 
is a Kleshchev multipartition. 
Now we can compute the decomposition matrix of 
$\ol{\CB}$ 
by the Jantzen sum formula, 
and its matrix is  
completely the same as the decomposition matrix of 
$\ol{\CB}$ 
in Proposition \ref{infinite S3}, 
replacing 
$k$ by $k+1$.  
Thus, 
by 
Lemma \ref{lem infinite cellular}, 
$\ol{\CB}$ is of infinite type,  
thus 
$\CB$ 
and so 
$\He_{n,4}(q,Q_1,\cdots,Q_4)$ is of infinite type.

For the case where $n \leq 3$, 
one can check case by case 
as in the proof of Proposition \ref{infinite S3}.   
\end{proof}

Now, we can prove Theorem \ref{main thm}. 

\begin{proof}[(Proof of Theorem \ref{main thm})]
The \lq\lq\, if " part 
is already shown in Proposition \ref{finiteness}. 
We show the \lq\lq only if " part, 
and 
we prove only the statement for $\He_{n,r}$ 
since the statement for $\Sc_{n,r}$ follows from one for $\He_{n,r}$ 
by Lemma \ref{lem Hecke Schur}. 

By Remarks \ref{remark} (\roiii), 
we may assume that 
$0 \leq f_1 \leq f_2 \leq \cdots \leq f_r \leq e-1$.

First, we consider the case where  
$r=2$. 
If 
$n \geq \min\{e, 2f^{+1}(Q_1,Q_2)+4\}$, 
then 
$\He_{n,2}(q,Q_1,Q_2)$ 
is of infinite type by 
\cite[Theorem 1.4]{AM04}.

Next, we consider the case where  
$r\geq 3$. 
If 
$n \geq \min \{e, 2f^{+1}(Q_1,\cdots,Q_r)+4\}$, 
then 
we have 
$n \geq \min \{e, 2(f_{i+1}-f_i)+4\}$ 
for some $i=1,\cdots,r-1$ 
or 
$n \geq \min \{e, 2(e-f_r)+4\}$.  
Take such $i$ 
(put $i=r$ if the last case occurs), 
then 
$\He_{n,2}(q,Q_i,Q_{i+1})$ (put $i+1=1$ if $i=r$) 
has infinite type by the above arguments. 
Thus, 
$\He_{n,r}(q,Q_1,\cdots,Q_r)$ 
has infinite type 
by Corollary \ref{cor induction}. 

If 
$n < \min\{e,2f^{+1}(Q_1,\cdots,Q_r)+4\}$ 
and $n \geq f^{+2}(Q_1,\cdots,Q_r)+1$, 
then 
there exist $i \in \{1,\cdots,r\}$ 
such that 
$\He_{n,3}(q,Q_i,Q_{i+1},Q_{i+2})$ 
(put $r+1=1,r+2=2$ if $i=r-1$ or $i=r$) 
satisfies 
the assumption in the Proposition \ref{infinite S3} 
(by adjusting the parameters by 
a scalar multiplication if necessary). 
By Proposition \ref{infinite S3}, 
such 
$\He_{n,3}(q,Q_i,Q_{i+1},Q_{i+2})$ 
has infinite type. 
Thus, 
we see that  
$\He_{n,r}(q,Q_1,\cdots,Q_r)$ 
has infinite type 
by Corollary \ref{cor induction}. 

Finally, 
we consider the case where  
$n<\min \big\{ 2f^{+1}(Q_1,\cdots,Q_r)+4,\,f^{+2}(Q_1,\cdots,\\  Q_r)+1\big\}$ 
and 
$n \geq g(Q_1,\cdots,Q_r)+2$. 
If 
$g(Q_1,\cdots,Q_r)=g_i+g_{i+1}$ 
for some 
$i\in \{1,\cdots,r\}$, 
we have 
$g(Q_1,\cdots,Q_r)=f^{+2}(Q_1,\cdots,Q_r)$. 
In this case, 
the conditions 
$n<f^{+2}(Q_1,\cdots,Q_r)+1$ 
and 
$n \geq g(Q_1,\cdots,Q_r)+2$
are not compatible.  
Thus, there exist 
$i, j \in \{1,\cdots,r\}$ 
such that 
$j-i>1$, 
$\{i,j\}\not=\{1,r\}$,  
and 
that 
$\He_{n,4}(q,Q_i,Q_{i+1},Q_{j},Q_{j+1})$ 
(put $r+1=1$ if $j=r$) 
satisfies the assumption in  
Proposition \ref{infinite condition g} 
(by adjusting the parameters by a scalar multiplication if necessary). 
By 
Proposition \ref{infinite condition g}, 
such 
$\He_{n,4}(q,Q_i,Q_{i+1},Q_{j},Q_{j+1})$ 
has infinite type. 
Thus 
$\He_{n,r}(q,Q_1,\cdots,Q_r)$ 
is  of infinite type 
by Corollary \ref{cor induction}.

This proves the \lq\lq only if " part. 
We have completed the proof of Theorem \ref{main thm}. 
\end{proof}


\remarks\ 
\label{final remarks}

(\roi) 
The Poincar\'e polynomial of the complex reflection group 
$W=\FS_n \ltimes (\ZZ/ r\ZZ)^n$ 
is given as 
\[ P_W(t)= \prod_{i=1}^n \frac{t^{ir} -1}{t-1}.\] 
Now, 
we consider the Ariki-Koike algebra $\He_{n,r}$ with one parameter, 
namely, in the case where 
the first relation in the definition of Ariki-Koike algebra (see \ref{definition of Ariki-Koike}) is 
\begin{align}
\label{one parameter }
(T_0 -q)(T_0-\zeta )(T_0-\zeta^2) \cdots (T_0-\zeta^{r-1})=0,
\end{align}
where $\zeta$ is a primitive $r$-th root of unity. 
We assume that 
$q$ is a primitive $e$-th root of unity and that $r$ divides $e$. 
Then we have $\zeta=q^{\frac{e}{r}}$. 
In order to apply Theorem \ref{main thm}, 
we rewrite the relation (\ref{one parameter }) 
by changing the generator 
$T_0$ 
by 
$q^{-1} T_0$ 
as follows. 
\[(T_0-1)(T_0-q^{\frac{e}{r}-1})(T_0 - q^{\frac{2e}{r} -1}) \cdots (T_0 - q^{\frac{(r-1)e}{r}-1}).\]
In this case, 
the condition 
$n < \min\{ 2f^{+1}(Q_1 , \cdots,Q_r)+4, f^{+2}(Q_1,\cdots, Q_r)+1, g(Q_1, \\   \cdots,Q_r)+2 \}$ 
is equivalent to the condition 
$ n \leq \frac{2e}{r}$. 
Thus, by Theorem \ref{main thm}, 
the Ariki-Koike algebra $\He_{n,r}$ is of finite type 
if and only if 
$ n \leq \frac{2e}{r}$. 
Moreover, if $\He_{n,r}$ is not semisimple then we have $\frac{e}{r} \leq n$. 
The condition 
$\frac{e}{r} \leq n \leq \frac{2e}{r}$ 
is equivalent to 
the condition that 
$q$ (a primitive $e$-th root of unity) 
is a simple root of the Poincar\'e polynomial $P_W(t)$. 
This result is compatible with a generalization of Uno's conjecture for Hecke algebras (\cite{Uno92}). 

(\roii) 
One can check that, 
if $\Sc_{n,r}(q,Q_1,\cdots,Q_r)$ 
is of finite type, 
then 
the weight of each block of $\Sc_{n,r}(q,Q_1,\cdots,Q_r)$ 
(in the sense of Fayers \cite{Fay06}) 
is less than or equal to one 
by  
\cite[Proposition 3.5]{Fay06} 
and the definition of the weight of a block 
\cite[2.1]{Fay06} 
combined with 
Lemma \ref{Lemma fin}.
On the other hand,  
if the weight of a block of $\Sc_{n,r}$ 
is $0$, 
then such a block is semisimple 
by \cite[Theorem 4.1]{Fay06}. 
If the  weight of a block of $\Sc_{n,r}$ 
is one, 
then such a block is of finite type 
by \cite[Theorem 4.12]{Fay06} 
combined with Proposition \ref{decom matrix}. 
(These facts give an alternate proof of Proposition \ref{finiteness}.)
Hence, 
it is likely that 
a block of cyclotomic $q$-Schur algebra $\Sc_{n,r}$ 
is of (non-semisimple) finite type 
if and only if 
the weight of the block is equal to one.


\end{document}